\documentclass{amsart}

\usepackage{amssymb}
\usepackage{amsmath}

\title[Nonexistence of doubly periodic solutions]
{Nonexistence of small doubly periodic solutions for dispersive equations}
\author{David M. Ambrose}
\address{Department of Mathematics, Drexel University, 3141 Chestnut St., Philadelphia, PA 19104, 
USA}
\thanks{DMA gratefully acknowledges support from the National Science Foundation through grants
DMS-1008387 and DMS-1016267.}
\author{J. Douglas Wright}
\address{Department of Mathematics, Drexel University, 3141 Chestnut St., Philadelphia, PA 19104, 
USA}
\thanks{JDW gratefully acknowledges support from the National Science Foundation through
grant DMS-1105635.}

\newtheorem{theorem}{Theorem}
\newtheorem{lemma}[theorem]{Lemma}
\newtheorem{cor}[theorem]{Corollary}
\newtheorem{remark}[theorem]{Remark}
\newtheorem{proposition}[theorem]{Proposition}

\begin{document}

\newcommand{\primed}{\!\!\!{}^{\prime}}

\begin{abstract}We study the question of existence of time-periodic, spatially periodic solutions
for dispersive evolution equations, and in particular, we introduce a framework for demonstrating
the nonexistence of such solutions.  We formulate the problem so that
doubly periodic solutions correspond to fixed points of a certain operator.  We prove that this
operator is locally contracting, for almost every temporal period, 
if the Duhamel integral associated to the evolution exhibits a weak smoothing property.
This implies the nonexistence of nontrivial, small-amplitude time-periodic solutions for almost
every period if the smoothing property holds.  This can be viewed as a partial analogue of
scattering for dispersive equations on periodic intervals, since scattering in free space implies the
nonexistence of small coherent structures.  We use a normal form to
demonstrate the smoothing property on 
specific examples, so that it can be seen that there are indeed equations for which the hypotheses
of the general theorem hold.  The nonexistence result is thus
established through the novel combination of
small divisor estimates and dispersive smoothing estimates.
The examples treated include the Korteweg-de Vries equation and the Kawahara equation.
\end{abstract}

\maketitle

\section{Introduction}

In the absence of the ability to ``explicitly" compute solutions of the Cauchy problem for a nonlinear dispersive system
by some specialized technique particular to the equation at hand  (such as complete integrability), 
coherent structures often 
form the backbone for both qualitative and quantitative descriptions of the dynamics of the system.
Such structures, be they traveling waves, self-similar solutions, time-periodic solutions or some other 
sort of solution, give
great insight into the short-time behavior of the system and often provide possible states
towards which solutions trend as time goes to infinity. 

For dispersive equations in free space, many authors have proved scattering results; we cannot
hope to list all such results here, but a sampling is  \cite{christWeinstein}, \cite{ginibreVelo},
\cite{liu-scatter},  \cite{ponceVega}, \cite{strauss-dispersion}.  In \cite{strauss-dispersion}, for instance,
Strauss showed that for a generalized Korteweg-de Vries (KdV) equation and for a nonlinear
Schr\"{o}dinger equation, all sufficiently small solutions decay to zero.  There is no generally 
agreed upon meaning for scattering on periodic intervals, and one cannot expect decay of solutions.
However, decay of solutions implies the nonexistence of small-amplitude coherent
structures, and the nonexistence of small coherent structures is a question which can be studied
on periodic intervals.  In the present work, we study the nonexistence of small time-periodic solutions
for dispersive equations on periodic intervals.  We prove a general theorem, showing that the
existence of sufficiently strong dispersive smoothing effects implies the nonexistence of small
doubly periodic solutions for almost every temporal period; we then demonstrate the required 
dispersive smoothing for particular examples.  We have partially described the results of the present
work briefly in the announcement \cite{ambroseWright-CR}.

One method of constructing doubly periodic solutions of dispersive equations is to use Nash-Moser-type 
methods.  These methods typically work in a function space with periodic boundary conditions
in both space and time, so that the Fourier transform of the evolution equation can be taken 
in both variables.  Then, a solution is sought nearby to an equilibrium solution.  The implicit function
theorem cannot be used due to the presence of small divisors.  Instead, the small divisors are
compensated for by the fast convergence afforded by Newton's method.  A version of such 
arguments is now known as the Craig-Wayne-Bourgain method, after Craig, Wayne, and Bourgain
used such arguments to demonstrate existence of doubly periodic solutions for a number of 
equations, such as nonlinear wave and nonlinear Schr\"{o}dinger equations 
\cite{bourgain1}, \cite{bourgain2}, \cite{craigWayne}, \cite{wayne}.  

Such methods have since been extended to other systems, such as the irrotational gravity water
wave, on finite or infinite depth, by Plotnikov, Toland, and Iooss \cite{plotnikovToland},
\cite{ioossPlotnikovToland}, or irrotational gravity-capillary waves by Alazard and Baldi
\cite{alazardBaldi}.  Also, Baldi used such methods to demonstrate existence of doubly periodic
solutions for perturbations of the Benjamin-Ono equation \cite{baldi}.

The typical result of these small divisor methods is the existence of small-amplitude doubly
periodic waves for the system under consideration, for certain values of the relevant parameters.
One such parameter is the frequency (or equivalently, the temporal period) of the solution; other
parameters may arise in specific applications, such as the surface tension parameter in 
\cite{alazardBaldi}.  With these methods, the parameter values 
for which solutions are shown to exist are typically in a Cantor set.

For certain completely integrable equations, time-periodic waves can be shown to exist by producing 
exact, closed-form solutions.  This is the case, for instance, for the KdV equation 
\cite{dubrovin} and the Benjamin-Ono equation \cite{dobrokhotovKrichever}, \cite{matsuno}, 
\cite{satsumaIshimori}.  The first author and Wilkening found that these
time-periodic solutions of the Benjamin-Ono equation form continuous families 
\cite{ambroseWilkening1}, \cite{ambroseWilkening2}, \cite{wilkening}; this is in sharp contrast
to results proved by small divisor methods, then.  The small divisor methods 
do not address the question of whether such results are optimal; it is not possible at present to conclude
whether or not time-periodic solutions exist as continuous families for such equations.

Thus, we find it of significant interest to develop further tools to answer the questions of both existence
and nonexistence of doubly periodic waves for dispersive partial differential equations.  In the present
contribution, we develop a framework by which the nonexistence of small-amplitude time-periodic 
waves can be established.  We do this by first using the Duhamel formula together with the time-periodic
ansatz.  Then we rewrite this formula, factoring out the linear operator.  The resulting equation yields 
the notion of time-periodic solutions as fixed points of a new operator, which is the composition
of the inverse of the linear operator with the Duhamel integral.

We then prove a general theorem, showing that under certain conditions, this operator
is contracting in a neighborhood of the origin in a certain function space.  Since the operator is the
composition of two operators, we prove estimates for these individual operators separately.  We are
able to prove that the inverse of the linear operator acts like differentiation of order $p$ for some
$p>1.$  Of course, in order to have the contracting property, the composition must map from 
some function space, $X,$ to itself.  We use Sobolev spaces, so since the inverse of the linear 
operator acts like differentiation of order $p,$ the Duhamel integral must satisfy an estimate
with some form of a gain at least $p$ 
derivatives.  In Theorem \ref{generalTheorem}, then, we have a general condition for the
nonexistence of small-amplitude time-periodic waves for almost every possible temporal period: if
the Duhamel integral possesses a weak form of smoothing (with an associated estimate), then
the equation does not possesss arbitrarily small doubly periodic waves for almost every period.
The results described thus far are the content of Section \ref{generalSection}.
We mention that we are aware of one other result in the literature on nonexistence of small
time-periodic solutions for almost every period; this is the paper \cite{delaLlaveWave}, in which
de la Llave uses a variational method to demonstrate nonexistence of small doubly periodic solutions
for nonlinear wave equations.

Clearly, we must address the question of whether there is an equation for which the truth of the
hypotheses of Theorem \ref{generalTheorem} can be demonstrated.
We demonstrate the required smoothing for the Duhamel integral associated to 
some dispersive evolution equations with fifth-order dispersion in Sections \ref{section:fifthOrder}
and \ref{section:kawahara}, and with seventh-order dispersion in Section \ref{section:seventhOrder}.
For these equations, the Duhamel integral satisfies a strong smoothing property: the 
Duhamel integral gains $p>1$ derivatives, compensating for the loss of derivatives from the
inverse of the linear operator.  Theorem \ref{generalTheorem} does not, however, require so strong
a smoothing property.  In Section \ref{section:kdv}, we demonstrate that the Duhamel integral for the
KdV equation satisfies a weaker smoothing property, allowing Theorem \ref{generalTheorem}
to be applied and demonstrating the nonexistence of small doubly periodic solutions for almost
every temporal period.

The estimates for the inverse of the linear operator (Lemma 
\ref{homogeneousLemma},
Lemma \ref{linearLemma2}, and Corollary \ref{almostEveryCorollary} below), in which we demonstrate
that the inverse of our linear operator acts like differentiation of order $p,$ are proved by small-divisor
techniques.  In fact, these are versions of classical results, such as can be found, for instance, 
in \cite{ghys}. 
As with all small-divisor results, some parameter values are discarded; in the present
case, the parameter is the temporal period of the solution.  Thus, we arrive at a result about 
nonexistence of small solutions for almost every possible temporal period.
Even though such a small divisor argument is classical, we provide our own proof because the 
detailed information about the set of temporal periods in the proof is helpful.

For our particular examples of dispersive equations, we prove smoothing estimates for the
Duhamel integral by following the lines of an argument by Erdogan and Tzirakis \cite{erdoganTzirakis}.
By using a normal form representation, Erdogan and Tzirakis showed that the 
Duhamel integral for the KdV equation gains $1-\varepsilon$ derivatives as compared to the initial
data, for any $\varepsilon>0.$  (We mention that a similar result has been demonstrated on the 
real line by Linares and Scialom \cite{linaresScialom}.)
For the equations with fifth-order dispersion which we consider in Sections \ref{section:fifthOrder}
and \ref{section:kawahara}, we find a gain of two derivatives.  This is in line with the usual,
expected gain of regularity from dispersion.  In general, if the dispersion relation is of order $r$
(say, for instance, the linearized evolution equation is, in the Fourier transform, 
$\hat{w}_{t}=ik^{r}\hat{w}$), then one expects to gain $(r-1)/2$ derivatives in some sense 
\cite{KPV-generalSmoothing}; this is known as the Kato smoothing effect \cite{kato}. 
With fifth-order dispersion, this means the expected gain is two derivatives.
Given the results of \cite{erdoganTzirakis} as well as the present work, it does appear
that it is reasonable to expect the same order of smoothing in the spatially periodic
setting, for the Duhamel integral.  In fact, in Section \ref{section:seventhOrder}, we show that
for an evolution equation with seventh-order dispersion, the gain of regularity on the Duhamel integral is 
four derivatives; thus, the smoothing effect in \cite{erdoganTzirakis} and the present work is not
the same as Kato smoothing, but is still due to the presence of dispersion.
To demonstrate our weaker smoothing estimate for the KdV equation in Section \ref{section:kdv},
we begin with the same normal form as before, but we estimate the terms differently.

We close with some discussion in Section \ref{section:discussion}.

\section{Nonexistence of doubly periodic solutions}\label{generalSection}

We begin with the evolution equation 
\begin{equation}\label{genericEvolution}
\partial_{t}u=Au+Nu,
\end{equation}
where $A$ is a linear operator and $N$ is a nonlinear operator.
Then, the solution of (\ref{genericEvolution}) with initial data $u(\cdot,t)=u_{0},$ 
if there is a solution, can be represented with the
usual Duhamel formula,
\begin{equation}\label{duhamel1}
u(\cdot,t)=e^{At}u_{0}+\int_{0}^{t}e^{A(t-\tau)}N(u(\cdot,\tau))\ d\tau.
\end{equation}
Given a time $t,$ we define the linear solution operator $S_{L}(t)=e^{At}$ and the difference
of the solution operator and the linear solution operator to be $S_{D}(t);$ thus, $S_{D}$ is
exactly the Duhamel integral:
\begin{equation}\nonumber
S_{D}(t)u_{0}=\int_{0}^{t}e^{A(t-\tau)}N(u(\cdot,\tau))\ d\tau.
\end{equation}

We work in the spatially periodic case, so we assume that solutions $u$ of (\ref{genericEvolution})
satisfy $$u(x+2\pi,t)=u(x,t), \qquad \forall x\in\mathbb{R}.$$
We assume that (\ref{genericEvolution}) maintains the mean of solutions; that is
given any $u$ in a reasonable function space, we have $$\int_{0}^{2\pi}u(x,t)\ dx =
\int_{0}^{2\pi}u_{0}(x)\ dx.$$  For the remainder of the present section, we will assume
that $u_{0}$ (and thus $u(\cdot,t)$) has mean zero.  

If $u_{0}$ is the initial data for a time-periodic solution of (\ref{genericEvolution}) with temporal period
$T,$ then we have
\begin{equation}\nonumber
u_{0}=S_{L}(T)u_{0}+S_{D}(T)u_{0}.
\end{equation}
We rewrite this as 
\begin{equation}\label{duhamel4}
(I-S_{L}(T)-S_{D}(T))u_{0}=0.
\end{equation}

Our goal is to demonstrate nonexistence of small-amplitude doubly periodic solutions of 
(\ref{duhamel4}), for certain temporal periods.  We begin now to focus only on certain values of $T,$ 
and our first restriction on values of $T$ is to ensure that $I-S_{L}(T)$ is invertible.  For $s>0,$
define the space $H^{s}_{0}$ to be the subset of the usual spatially periodic 
Sobolev space $H^{s},$ such that for
all $f\in H^{s}_{0},$ the mean of $f$ is equal to zero.  We assume that the operator $S_{L}(t)$ is 
bounded, $$S_{L}(t):H^{s}_{0}\rightarrow H^{s}_{0},\qquad \forall t\in\mathbb{R}.$$
Then, we define the set $W$ to be
$$W=\{t\in(0,\infty): \mathrm{ker}(I-S_{L}(t))=\{0\}\}.$$  For any $T\in W,$ we rewrite (\ref{duhamel4})
by factoring out $I-S_{L}(T):$
\begin{equation}\nonumber
(I-S_{L}(T))(I-(I-S_{L}(T))^{-1}S_{D}(T))u_{0}=0.
\end{equation}

We see, then, that if $u_{0}$ is the initial data for a time-periodic solution of (\ref{genericEvolution})
with temporal period $T\in W,$ then $u_{0}$ is a fixed point of the operator
$$K(T):=(I-S_{L}(T))^{-1}S_{D}(T).$$
If we can show that this is (locally) a contraction on $H^{s}_{0},$ then there are no (small) nontrivial
time-periodic solutions in the space $H^{s}_{0}$ with temporal period $T.$  To establish this,
we will need estimates both for $(I-S_{L}(T))^{-1}$ and for $S_{D}(T).$  In Sections \ref{section:linear1}
and \ref{section:linear2}, we establish estimates for $(I-S_{L}(T))^{-1};$ the results are that the
symbol can be bounded as $|k|^{p},$ where $k$ is the variable in Fourier space, for some $p>1,$
for certain values of $T.$
Thus, the inverse of the linear operator acts like differentiation of order $p>1.$  In Section
\ref{section:generalTheorem}, then, we will state a corollary of these estimates: if the operator 
$S_{D}$ satisfies a certain estimate related to a gain of $p$ derivatives, then 
the operator $K(T)$ is locally contracting, and thus there are no small time-periodic solutions
with temporal period $T.$

\subsection{The linear estimate: the homogeneous case}\label{section:linear1}

In this section, we prove our estimate for $(I-S_{L}(T))^{-1}$ in the case that the linear operator $A$
has symbol 
\begin{equation}\label{homogeneousSymbol}
\mathcal{F}(A)(k)=ik^{r}, \qquad \forall k\in\mathbb{Z}.
\end{equation}
This estimate is the content of
Lemma \ref{homogeneousLemma}.  We note that this is not, strictly speaking, useful to us, 
as we will prove a version of the lemma for
a more general class of operators $A$ in the following section.  However, we include Lemma
\ref{homogeneousLemma} because the simplicity of the form of (\ref{homogeneousSymbol}) allows
for clarity of exposition.  The proof of Lemma \ref{linearLemma2} below will build upon the proof of
Lemma \ref{homogeneousLemma}, which we now present.

\begin{lemma}\label{homogeneousLemma}
 Let the linear operator $A$ be given by (\ref{homogeneousSymbol}).  
Let $0<T_{1}<T_{2}$ be given.  Let $0<\delta<T_{2}-T_{1}$ be given.  Let $p>1$ be given.
There exists a set $W_{p,\delta}\subseteq[T_{1},T_{2}]\cap W$ 
and there exists $c_{1}>0$ such that the Lebesgue
measure of $W_{p,\delta}$ satisfies $\mu(W_{p,\delta})>T_{2}-T_{1}-\delta,$ and for all 
$k\in\mathbb{Z}\setminus\{0\},$ for all $T\in W_{p,\delta},$ we have
$$\left|\mathcal{F}(I-S_{L}(T))^{-1}(k)\right|<c_{1}|k|^{p}.$$
\end{lemma}

\noindent{\bf{Proof:}}  For the moment, we fix $k\in\mathbb{Z}\setminus\{0\}.$
We need to estimate 
\begin{equation}\label{toEstimate-Linear}
\left|\mathcal{F}(I-S_{L}(T))^{-1}(k)\right|=
\left|\frac{1}{1-\exp\{ik^{r}T\}}\right|=\frac{1}{\sqrt{2}}(1-\cos(k^{r}T))^{-1/2}.
\end{equation}
Clearly, the symbol of the inverse operator
is undefined if there exists $n\in\mathbb{N}$ such that $T=\frac{2\pi n}{|k|^{r}}.$  With the assumption 
that $T\in[T_{1},T_{2}],$ the associated values of $n$ comprise the set
$\mathcal{N}:=\left[\frac{|k|^{r}T_{1}}{2\pi},\frac{|k|^{r}T_{2}}{2\pi}\right]\cap\mathbb{N}.$  We remove a small set of possible periods around 
each of these
values; that is to say, we consider
\begin{equation}\label{setOfTimes}
T\in[T_{1},T_{2}]\setminus\bigcup_{n\in\mathcal{N}}\left[
\frac{2\pi n}{|k|^{r}}-\varepsilon, \frac{2\pi n}{|k|^{r}}+\varepsilon\right],
\end{equation}
for some $0< \varepsilon\ll 1$ to be specified.

To start, we may notice that the collection of intervals $\left[\frac{2\pi n}{|k|^{r}}-\varepsilon,
\frac{2\pi n}{|k|^{r}}+\varepsilon\right]$ do not overlap for different values of $n$ as long as, 
for all $n\in\mathcal{N},$ we have 
$$\frac{2\pi n}{|k|^{r}}+\varepsilon<\frac{2\pi(n+1)}{|k|^{r}}-\varepsilon.$$
This is satisfied as long as $\varepsilon<\frac{\pi}{|k|^{r}}.$  When we choose $\varepsilon,$ this
condition will be satisfied.  

A simple calculation shows that, on an interval of values of $\theta$ which does not include
an integer multiple of $2\pi,$ the value of $(1-\cos(\theta))$ is minimized at the endpoints of the
interval.  Thus, the value of (\ref{toEstimate-Linear}) is largest on our set of possible periods
at the values $T=\frac{2\pi n}{|k|^{r}}\pm\varepsilon,$ for $n\in\mathcal{N}.$  At such values, we
find
\begin{multline}\label{cosineToBeEstimated}
\left|\mathcal{F}\left(\left(I-S_{L}\left(\frac{2\pi n}{|k|^{r}}\pm\varepsilon\right)\right)^{-1}\right)(k)\right|
=\frac{1}
{\sqrt{2}\left(1-\cos\left(k^{r}\left(\frac{2\pi n}{|k|^{r}}\pm\varepsilon\right)\right)\right)^{1/2}}
\\
=\frac{1}{\sqrt{2}(1-\cos(\pm2\pi n \pm \varepsilon k^{r}))^{1/2}}
=\frac{1}{\sqrt{2}(1-\cos(\varepsilon k^{r}))^{1/2}}.
\end{multline}
We will now perform a Taylor expansion for cosine, paying attention to the error estimates.
For any $\theta\in\mathbb{R},$ 
we have the formula $$\cos(\theta)=1-\frac{\theta^{2}}{2}+\frac{\sin(\xi)\theta^{3}}{6},$$ for
some $\xi$ between $0$ and $\theta.$
We notice that $\left|\frac{\sin(\xi)\theta^{3}}{6}\right|\leq \frac{\theta^{2}}{4},$ as long as 
$|\theta|\leq\frac{3}{2}.$  Thus, for $|\theta|\leq\frac{3}{2},$ we have
\begin{equation}\label{cosineEstimate}
\left|1-\cos(\theta)\right|=\left|\frac{\theta^{2}}{2}-\frac{\sin(\xi)\theta^{3}}{6}\right|
\geq\frac{\theta^{2}}{4}.
\end{equation}
Combining (\ref{cosineEstimate}) with (\ref{cosineToBeEstimated}), we find that for any 
$T$ satisfying (\ref{setOfTimes}), we have
\begin{equation}\label{almostFinalLinear1}
\left|\mathcal{F}\left((I-S_{L}(T))^{-1}(k)\right)\right|\leq \frac{1}{(\sqrt{2})(\frac{\varepsilon |k|^{r}}{2})}
=\frac{\sqrt{2}}{\varepsilon |k|^{r}}.
\end{equation}

We choose $\varepsilon=c_{0}|k|^{-p-r}.$ This choice of $\varepsilon$ immediately yields
the claimed estimate for the symbol.  Recall that we have specified $p>1;$ the positive 
constant $c_{0}$ is to be specified.
The conditions we have placed on $\varepsilon$ above are (1) $\varepsilon<\frac{\pi}{|k|^{r}}$ and
(2) $\varepsilon|k|^{r}<\frac{3}{2}.$  These conditions are both 
satisfied as long as $c_{0}<\frac{3}{2}.$

For fixed $k\in\mathbb{Z}\setminus\{0\}$ 
and for fixed $n\in\mathcal{N},$ we have removed a set of measure
$2\varepsilon=2c_{0}|k|^{-r-p}$ from the interval $[T_{1},T_{2}].$  
Since $\mathcal{N}$ is the intersection of
an interval with the natural numbers, we see that the cardinality of $\mathcal{N}$ less than
$|k|^{r}(T_{2}-T_{1})+1,$ so for fixed $k,$ we have removed a set of measure no more than
$2c_{0}(1+T_{2}-T_{1})|k|^{-p}.$  
Summing over $k,$ since we have chosen $p>1,$ the measure of the set
we have removed is finite.  Taking $c_{0}$ sufficiently small, we can conclude that the measure
of the set which is removed has Lebesgue measure smaller than $\delta.$  To be definite, we
write the definition of the set $W_{p,\delta},$ which is 
$$W_{p,\delta}=[T_{1},T_{2}] \setminus
\bigcup_{k\in\mathbb{Z}\setminus\{0\}}\bigcup_{n\in\mathcal{N}}
\left[\frac{2\pi n}{|k|^{r}}-\frac{c_{0}}{|k|^{r+p}}, \frac{2\pi n}{|k|^{r}}+\frac{c_{0}}{|k|^{r+p}}\right],$$
where $c_{0}$ is chosen so that $0<c_{0}<\frac{3}{2},$ and also so that
$$c_{0}<\frac{\delta}{2(1+T_{2}-T_{1})\displaystyle\sum_{k\in\mathbb{Z}\setminus\{0\}}|k|^{-p}}.$$
Finally, the constant $c_{1}$ is given by $c_{1}=\displaystyle\frac{\sqrt{2}}{c_{0}}.$
This completes the proof. \hfill$\blacksquare$

We note that we can see clearly the dependence of the constant $c_{1}$ on the parameters
$T_{1},$ $T_{2},$ $p,$ and $\delta.$  If we want a larger set of potential periods, we could take
$T_{2}-T_{1}$ larger or $\delta$ smaller; this would result in a larger value of $c_{1}.$  Choosing smaller
values of $p>1$ also leads to larger values of $c_{1}.$

\subsection{The linear estimate: the nonresonant case}\label{section:linear2}
 The estimate of Section \ref{section:linear1}
can be generalized to allow operators which include lower-order terms, in some cases:
there must not be a resonance between the different terms, in the sense that we require
$\mathcal{F}A(k)\neq0,$ for any nonzero $k\in\mathbb{Z}.$ To be very precise, we consider
linear operators $A$ which satisfy the conditions {\bf(H)}, which we now describe.\\

\noindent{\bf(H)}
Let $M\in\mathbb{N}$ be given, with $M\geq 2.$  For all $m\in\{1,2,\ldots,M\},$ let $r_{m}\in\mathbb{R}$
such that $r_{1}>r_{2}>\ldots>r_{M}>0$ be given.
For all $m\in\{1,2,\ldots,M\},$ let
$Z_{m}\subseteq\mathbb{R}$ be bounded.  Let $Z=Z_{1}\times Z_{2}\times\ldots Z_{M}.$
Let $\vec{\alpha}=(\alpha_{1},\alpha_{2},\ldots,\alpha_{M}).$
Assume there exists $\beta_{1}>0$ and $\beta_{2}>0$ such that for all $\alpha_{1}\in Z_{1},$
\begin{equation}\label{h1}
|\alpha_{1}|>\beta_{1},
\end{equation} 
and for all $\vec{\alpha}\in Z,$ for all $k\in\mathbb{Z}\setminus\{0\},$
\begin{equation}\label{h2}
\left|\sum_{m=1}^{M}\alpha_{m}k^{r_{m}}\right|\geq \beta_{2}.
\end{equation}
Given $\vec{\alpha}\in Z,$ let the linear operator $A$ be defined 
through its symbol as 
\begin{equation}\label{nonresonantA}
\mathcal{F}(A)(k)=i\sum_{m=1}^{M}\alpha_{m}k^{r_{m}},\qquad\forall k\in\mathbb{Z}.
\end{equation}

We remark that the condition (\ref{h1}) ensures that the equation is dispersive of order $r_{1}$ (i.e., the leading-order
term is of the same order for all $\vec{\alpha}\in Z$).  The condition (\ref{h2}) ensures that the symbol never
vanishes.

\begin{lemma}\label{linearLemma2}
Let the set $Z$ and the linear operator $A$ satisfy the hypotheses {\bf(H)}, so that in particular, 
$A$ is defined by (\ref{nonresonantA}).
Let $0<T_{1}<T_{2}$ be given.  Let $0<\delta<T_{2}-T_{1}$ be given.  Let $p>1$ be given.
There exists a set $W_{p,\delta}\subseteq[T_{1},T_{2}]\cap W$ 
and there exists $c_{1}>0$ such that the Lebesgue
measure of $W_{p,\delta}$ satisfies $\mu(W_{p,\delta})>T_{2}-T_{1}-\delta,$ and 
such that for all $\vec{\alpha}\in Z,$ for all 
$k\in\mathbb{Z}\setminus\{0\},$ for all $T\in W_{p,\delta},$ we have
$$\left|\mathcal{F}(I-S_{L}(T))^{-1}(k)\right|<c_{1}|k|^{p}.$$
\end{lemma}

\noindent{\bf{Proof:}}
The proof of Lemma \ref{homogeneousLemma} can be repeated, with $k^{r}$
replaced in every instance by $\displaystyle\sum_{m=1}^{M}\alpha_{m}k^{r_{m}},$
until (\ref{almostFinalLinear1}).  We then choose $\varepsilon=c_{0}|k|^{-p-r_{1}},$ with $c_{0}$ to
be specified.  

Similarly to the previous case, we have the conditions (i) 
$\varepsilon<(\pi)\left|\displaystyle\sum_{m=1}^{M}\alpha_{m}k^{r_{m}}\right|^{-1},$ and (ii)
$\varepsilon\left|\displaystyle\sum_{m=1}^{M}\alpha_{m}k^{r_{m}}\right|<\displaystyle\frac{3}{2}.$
The relevant product satisfies the following inequality:
$$\left|\varepsilon\sum_{m=1}^{M}\alpha_{m}k^{r_m}\right|=c_{0}\left|\sum_{m=1}^{M}\alpha_{m}
k^{r_{m}-r_{1}-p}\right|\leq c_{0}\sum_{m=1}^{M}|\alpha_{m}|,\qquad 
\forall k\in\mathbb{Z}\setminus\{0\}.$$
Thus, recalling hat the sets $Z_{m}$ are bounded,
we can clearly take $c_{0}$ sufficiently small to satisfy conditions (i) and (ii).

We must revisit the estimate (\ref{almostFinalLinear1}) in the present setting.  We have
\begin{equation}\nonumber
|\mathcal{F}(I-S_{L}(T))^{-1}(k)|\leq
\frac{\sqrt{2}}{\varepsilon\left|\sum_{m=1}^{M}\alpha_{m}k^{r_{m}}\right|}
=\left(\frac{\sqrt{2}}{c_{0}}\right)\frac{|k|^{p}}
{\left|\sum_{m=1}^{M}\alpha_{m}k^{r_{m}-r_{1}}\right|}.
\end{equation}
Since $\displaystyle\sum_{m=1}^{M}\alpha_{m}k^{r_{m}}$ is never equal to zero, we conclude that
$\displaystyle\sum_{m=1}^{M}\alpha_{m}k^{r_{m}-r_{1}}$ is also never equal to zero, for 
$k\in\mathbb{Z}\setminus\{0\}.$  
Furthermore, there exists $K\in\mathbb{N}$ such that for all $k\in\mathbb{Z}$ satisfying $|k|>K,$
for all $\vec{\alpha}\in Z,$ we have
$$\left|\sum_{m=2}^{M}\alpha_{m}k^{r_{m}-r_{1}}\right|\leq
|K|^{r_{2}-r_{1}}\sum_{m=2}^{M}|\alpha_{m}|
<\frac{\beta_{1}}{2}.$$
This implies that, for all $k\in\mathbb{Z}$ with $|k|>K,$ for all $\vec{\alpha}\in Z,$
$$\left|\sum_{m=1}^{M}\alpha_{m}k^{r_{m}-r_{1}}\right|=\left|\alpha_{1}
+\sum_{m=2}^{\infty}\alpha_{m}k^{r_{m}-r_{1}}\right|\geq \frac{\beta_{1}}{2}>0,$$
$$\left|\sum_{m=1}^{M}\alpha_{m}k^{r_{m}-r_{1}}\right|=\left|\alpha_{1}
+\sum_{m=2}^{\infty}\alpha_{m}k^{r_{m}-r_{1}}\right|\leq 
\left(\sup_{\alpha_{1}\in Z_{1}}|\alpha_{1}|\right)+\frac{\beta_{1}}{2}.$$
Furthermore, for any $k\in\mathbb{Z}\setminus\{0\}$ satisfying $|k|\leq K,$ we know
$$\inf_{\vec{\alpha}\in Z}\frac{|\sum_{m=1}^{M}\alpha_{m}k^{r_{m}}|}{|k^{r_{1}}|}
=\inf_{\vec{\alpha}\in Z}\left|\sum_{m=1}^{M}\alpha_{m}k^{r_{m}-r_{1}}\right|\in(0,\infty),$$
$$\sup_{\vec{\alpha}\in Z}\frac{|\sum_{m=1}^{M}\alpha_{m}k^{r_{m}}|}{|k^{r_{1}}|}
=\sup_{\vec{\alpha}\in Z}\left|\sum_{m=1}^{M}\alpha_{m}k^{r_{m}-r_{1}}\right|\in(0,\infty).$$
Combining this information, we see that $$\inf_{k\in\mathbb{Z}\setminus\{0\}} \inf_{\vec{\alpha}\in Z}
\left|\sum_{m=1}^{M}\alpha_{m}k^{r_{m}-r_{1}}\right|\in(0,\infty),$$
$$\sup_{k\in\mathbb{Z}\setminus\{0\}}\sup_{\vec{\alpha}\in Z}
\left|\sum_{m=1}^{M}\alpha_{m}k^{r_{m}-r_{1}}\right|\in(0,\infty).$$
Our value of $c_{1}$ is therefore
$$c_{1}=\left(\frac{\sqrt{2}}{c_{0}}\right)\left(
\inf_{k\in\mathbb{Z}\setminus\{0\}}\inf_{\vec{\alpha}\in Z}
\left|\sum_{m=1}^{M}\alpha_{m}k^{r_{m}-r_{1}}\right|\right)^{-1}.$$

In the current setting, for each $\vec{\alpha}\in Z,$ for each $k\in\mathbb{Z}\setminus\{0\},$
the set $\mathcal{N}$ is defined by
$$\mathcal{N}=\left[\frac{T_{1}\left|\sum_{m=1}^{M}\alpha_{m}k^{r_{m}}\right|}{2\pi},
\frac{T_{2}\left|\sum_{m=1}^{M}\alpha_{m}k^{r_m}\right|}{2\pi}\right]\cap\mathbb{N}.$$
Thus, for all $\vec{\alpha}\in Z,$ for all $k\in\mathbb{Z}\setminus\{0\},$
the cardinality of $\mathcal{N}$ satisfies
\begin{multline}\nonumber
\mathrm{card}(\mathcal{N})\leq(T_{2}-T_{1})\left|\sum_{m=1}^{M}\alpha_{m}k^{r_{m}}\right|+1
\leq (T_{2}-T_{1})|k|^{r_{1}}\left|\sum_{m=1}^{M}\alpha_{m}k^{r_{m}-r_{1}}\right|+1
\\
\leq (T_{2}-T_{1})|k|^{r_{1}}\left(\sup_{k\in\mathbb{Z}\setminus\{0\}}\sup_{\vec{\alpha}\in Z}
\left|\sum_{m=1}^{M}\alpha_{m}k^{r_{m}-r_{1}}\right|\right)+1.
\end{multline}
We then take the product $2\varepsilon\mathrm{card}(\mathcal{N}),$ finding the estimate
$$2\varepsilon\mathrm{card}(\mathcal{N})\leq 2c_{0}\left(
1+(T_{2}-T_{1})\sup_{k\in\mathbb{Z}\setminus\{0\}}\sup_{\vec{\alpha}\in Z}
\left|\sum_{m=1}^{M}\alpha_{m}k^{r_{m}-r_{1}}\right|
\right)|k|^{-p}.$$
As in the proof of Lemma \ref{homogeneousLemma}, we sum over $k,$ and we find that we can
take $c_{0}$ sufficiently small to satisfy the remaining claims of the lemma.
\hfill$\blacksquare$ 

Lemma \ref{linearLemma2} allows for two kinds of uniformity:
the same set $W_{p,\delta}$ works for all
$\vec{\alpha}\in Z,$ and the constant $c_{1}$ is able to be used
for all $T\in W_{p,\delta}.$  The cost of this uniformity with respect to the constant $c_{1}$
is that the set $W_{p,\delta}$ does not have 
full measure.  By sending $\delta$ to zero, we can achieve an estimate for almost every 
$T\in[T_{1},T_{2}],$ but then the constant will depend on
the choice of $T.$  In doing this, we are able to maintain the uniformity with respect to the set $Z.$
This is the content of the following corollary.

\begin{cor}\label{almostEveryCorollary}
Let the set $Z$ and the linear operator $A$ satisfy the hypotheses {\bf(H)}, with $A$ given by 
(\ref{nonresonantA}).
Let $p>1$ be given.  Let $0<T_{1}<T_{2}$ be given.
For almost every $T\in[T_{1},T_{2}],$ there exists $c>0$ such that for all $\vec{\alpha}\in Z$ and for
all $k\in\mathbb{Z}\setminus\{0\},$
we have the estimate 
$$\left|\mathcal{F}\left(\left(I-S_{L}(T)\right)^{-1}\right)(k)\right|\leq c|k|^{p}.$$
\end{cor}

\noindent{\bf{Proof:}} For any $\delta$ 
satisfying $0<\delta<T_{2}-T_{1},$ let the set $W_{p,\delta}$ be as in 
Lemma \ref{linearLemma2}.  For any $T\in[T_{1},T_{2}],$ if there exists a value of $\delta$ such that
$T\in W_{p,\delta},$ then the desired estimate is satisfied.  Let 
$$W_{p}=\bigcup_{n=2}^{\infty}W_{p,(T_{2}-T_{1})/n}.$$ Then, for all $T\in W_{p},$ the estimate holds.
Since for all $n\geq 2$ we have $W_{p,(T_{2}-T_{1})/n}\subseteq W_{p}\subseteq[T_{1},T_{2}],$ 
and since
$\mu(W_{p,(T_{2}-T_{1})/n})\geq (T_{2}-T_{1})(1-1/n),$ we see that 
$\mu(W_{p})=T_{2}-T_{1}.$  This completes the proof.
\hfill$\blacksquare$

\subsection{The general theorem}\label{section:generalTheorem}

We are now able to state a general theorem which follows from the above discussion.
In Theorem \ref{generalTheorem}, when we say that $N$ is ``as above,'' this includes the 
property that the evolution equation (\ref{genericEvolution}) preserves the mean value of the 
initial data.  The theorem contains two statements about nonexistence of small-amplitude 
time-periodic solutions.  The first statement is for a given $\delta>0;$ for $T$ in the set $W_{p,\delta},$
we conclude that there is a uniform threshold for the amplitude of time-periodic solutions.  For the
second statement, we conclude that for almost any $T\in[T_{1},T_{2}],$ there is a threshold for
the amplitude of time-periodic solutions; this second statement is not uniform.  These results
are conditional on the existence of smoothing estimates.  In Sections \ref{section:fifthOrder},
\ref{section:kawahara}, \ref{section:seventhOrder}, and \ref{section:kdv}, 
we will demonstrate the required smoothing estimate for particular equations.

\begin{theorem}\label{generalTheorem}
Let $0<T_{1}<T_{2}$ be given, and let $0<\delta<T_{2}-T_{1}$ be given. 
Let the set $Z$ and the operator $A$ satisfy the hypotheses {\bf(H)}, with
$A$ given by (\ref{nonresonantA}).
Let  the nonlinear operator $N$ be as above.  Assume there exists
$p>1,$ $\tilde{p}\geq0,$ $q>0,$ $s\geq0,$ $c>0,$ and $\eta>0$ such that for all $u_{0}\in H^{s+\tilde{p}}_{0}$ with
$\|u_{0}\|_{H^{s+\tilde{p}}}\leq \eta,$ for all $\vec{\alpha}\in Z,$ the following estimate is satisfied:
\begin{equation}\label{smoothingGeneral}
\|S_{D}(T)u_{0}\|_{H^{s+p}}\leq c\|u_{0}\|_{H^{s}}\|u_{0}\|_{H^{s+\tilde{p}}}^{q},
\end{equation}
for all $T\in W_{p,\delta}.$  Then, there exists $r_{0}>0$ such that if $u$ is
a smooth, nontrivial,
mean-zero time-periodic solution of (\ref{genericEvolution}) with temporal period $T\in W_{p,\delta},$
then 
\begin{equation}\label{goodConclusion}
\inf_{t\in[0,T]}\|u\|_{H^{s+\tilde{p}}}\geq r_{0}.\end{equation}
Furthermore, if (\ref{smoothingGeneral}) holds
for every $T\in W_{p},$ then for every $T\in W_{p},$ there exists 
$r_{0}>0$ such that if $u$ is
a nontrivial mean-zero time-periodic solution of (\ref{genericEvolution}) with temporal period $T,$
then (\ref{goodConclusion}) holds.
\end{theorem}

\noindent{\bf{Proof:}} The assumptions of the theorem, together with either Lemma \ref{linearLemma2}
or Corollary \ref{almostEveryCorollary}, imply that there exists $C>0$ such that for all $u_{0}\in H_{0}^{s+\tilde{p}}$ satisfying
$\|u_{0}\|_{H^{s+\tilde{p}}}<\eta,$
$$\|K(T)u_{0}\|_{H^{s}}\leq C\|u_{0}\|_{H^{s}}\|u_{0}\|_{H^{s+\tilde{p}}}^{q}.$$
If $u_{0}$ satisfies $$0<C\|u_{0}\|_{H^{s+\tilde{p}}}^{q}<1,$$
then $\|K(T)u_{0}\|_{H^{s}}<\|u_{0}\|_{H^{s}},$ and thus $u_{0}$ cannot be a fixed point of $K(T).$ 
Thus, the only fixed point in a ball around zero is zero.
\hfill$\blacksquare$

\begin{remark}In the announcement \cite{ambroseWright-CR}, the version of this theorem which
appeared was restricted to the case $\tilde{p}=0.$  In this case, the inverse of the linear operator
acts like differentiation of order $p>1,$ and the Duhamel integral has a compensating gain of $p$
derivatives.  In Sections \ref{section:fifthOrder} and \ref{section:further}, 
we will give examples for which this smoothing property
holds; these examples include the Kawahara equation.  
However, as the estimate (\ref{smoothingGeneral}) shows, what is needed is actually much
weaker than the Duhamel integral gaining $p$ derivatives; instead, it is only 
necessary that $s+p$ derivatives of the Duhamel integral satisfy a nonlinear estimate in which one 
factor involves only $s$ derivatives.  In Section \ref{section:kdv}, we will demonstrate that 
the estimate (\ref{smoothingGeneral}) holds with $\tilde{p}>0$ for the KdV equation.
\end{remark}

\section{Application to a fifth-order dispersive equation}\label{section:fifthOrder}

In this section, we will apply the above results to a specific dispersive equation, with sufficiently
strong dispersion, with $\tilde{p}=0.$  
We are using a version of the Erdogan-Tzirakis argument \cite{erdoganTzirakis}, 
for the equation
\begin{equation}\label{fifthOrderEquationTilde}
\partial_{t}\tilde{u}=\partial_{x}^{5}\tilde{u}-2\tilde{u}\partial_{x}\tilde{u}+\tilde{\omega}\partial_{x}\tilde{u},
\end{equation} (for any $\tilde{\omega}\in\mathbb{R}$)
 to get the desired smoothing effect.  As we have discussed, we consider the spatially periodic case,
 with spatial period equal to $2\pi.$  
 We first rewrite (\ref{fifthOrderEquationTilde}) to remove
 the mean, and also to remove the tildes.

We consider the initial condition $\tilde{u}(x,0)=\tilde{g}(x).$  Assume the mean of $\tilde{g}$ is equal to
$\bar{g},$ which can be any real number.  Let $u=\tilde{u}-\bar{g}.$  Since the evolution for
$\tilde{u}$ conserves the mean, the mean of $u$ will equal zero at all times.  The evolution
equation satisfied by $u$ is 
\begin{equation}\label{fifthOrderEquation}
\partial_{t}u=\partial_{x}^{5}u-2u\partial_{x}u+\omega\partial_{x}u,
\end{equation}
where $\omega=\tilde{\omega}-2\bar{g}.$  The initial data for (\ref{fifthOrderEquation}) is
$g=\tilde{g}-\bar{g},$ which of course has mean zero. 

We now discuss the appropriate existence theory for (\ref{fifthOrderEquation}).

\subsection{Existence Theory}\label{section:existenceTheory}

The well-posedness of the initial value problem for (\ref{fifthOrderEquation}) 
(or for (\ref{fifthOrderEquationTilde})) has been established in \cite{bourgain-wellposed}, in the
space $H^{s},$ for $s>0;$ a more general family of equations including these was also
shown to be well-posed in $H^{s}$ for $s>\frac{1}{2}$ in \cite{huLi}.  In the present work,
we are not concerned with demonstrating results at the lowest possible regularity, but instead
in finding estimates which will work with the nonexistence argument.  Towards this end, we
will give a simple existence theorem for the initial value problem for (\ref{fifthOrderEquation}) in
the space $H^{6},$ as the resulting estimates will be useful.  We mention that the choice of
$H^{6}$ as the function space is made so that the solutions are classical solutions.

\begin{proposition}
Let $u_{0}\in H^{6}$ be given.  There exists $T>0$ and a unique $u\in C([0,T];H^{6})$ such that 
$u$ solves the initial value problem (\ref{fifthOrderEquation}) with data $u(\cdot,0)=u_{0}.$
\end{proposition}

\noindent{\bf{Proof:}}
We begin by introducing a mollifier, $\mathcal{J}_{\varepsilon},$ for $\varepsilon>0.$
We use the mollifier to make an approximate evolution equation,
\begin{equation}\label{mollEvol}
\partial_{t}u^{\varepsilon}=\mathcal{J}_{\varepsilon}^{2}\partial_{x}^{5}u^{\varepsilon}-
\mathcal{J}_{\varepsilon}(2(\mathcal{J}_{\varepsilon}u^{\varepsilon})
(\mathcal{J}_{\varepsilon}u^{\varepsilon}_{x}))
+\omega\mathcal{J}_{\varepsilon}^{2}u^{\varepsilon}_{x}.\end{equation}
When combined with the mollifier, all of the derivatives on the right-hand side have
become bounded operators, and thus solutions for the initial value problem for $u^{\varepsilon},$
with initial data $u^{\varepsilon}(\cdot,0)=g\in H^{6},$ exist in $C([0,T_{\varepsilon}]; H^{6})$
by Picard's Theorem (cf. Chapter 3 of \cite{majdaBertozzi}).  

In order to show that the interval of existence can be taken to be independent of $\varepsilon,$
we must make an energy estimate.  We let the energy functional be an equivalent version of the
square of the $H^{6}$ norm:
$$\mathcal{E}(t)=\frac{1}{2}\int_{0}^{2\pi}(u^{\varepsilon})^{2}+(\partial_{x}^{6}u^{\varepsilon})^{2}\ dx.$$
Taking the time derivative, we find
$$\frac{d\mathcal{E}}{dt}=\int_{0}^{2\pi}(u^{\varepsilon})(\partial_{t}u^{\varepsilon})\ dx
+\int_{0}^{2\pi}(\partial_{x}^{6}u^{\varepsilon})(\partial_{t}\partial_{x}^{6}u^{\varepsilon})\ dx = I + II.$$

For $I,$ we plug in from the evolution equation, using the fact that $\mathcal{J}_{\varepsilon}$ is
self-adjoint and commutes with $\partial_{x}:$
\begin{equation}\nonumber
I=\int_{0}^{2\pi}(\mathcal{J}_{\varepsilon}u^{\varepsilon})
\partial_{x}^{5}(\mathcal{J}_{\varepsilon}u^{\varepsilon})
-2(\mathcal{J}_{\varepsilon}u^{\varepsilon})^{2}(\mathcal{J}_{\varepsilon}u^{\varepsilon}_{x})
+\omega(\mathcal{J}u^{\varepsilon})(\mathcal{J}_{\varepsilon}u^{\varepsilon}_{x})\ dx.
\end{equation}
All of these terms vanish upon integrating by parts and/or recognizing perfect derivatives; 
therefore, $I=0.$

To study the term $II,$ it is helpful to first apply six spatial derivatives to the equation (\ref{mollEvol}).
We use the product rule, finding the following:
\begin{multline}\label{withProductRule}
\partial_{t}\partial_{x}^{6}u^{\varepsilon}
=\mathcal{J}_{\varepsilon}^{2}\partial_{x}^{11}u^{\varepsilon}
-\mathcal{J}_{\varepsilon}(2(\mathcal{J}_{\varepsilon}u^{\varepsilon})
(\partial_{x}^{7}\mathcal{J}_{\varepsilon}u^{\varepsilon})+\omega\mathcal{J}_{\varepsilon}^{2}\partial_{x}^{7}
u^{\varepsilon}
\\
-2\mathcal{J}_{\varepsilon}\left(\sum_{m=1}^{6}{6\choose m}
(\partial_{x}^{m}\mathcal{J}_{\varepsilon}u^{\varepsilon})(
\partial_{x}^{7-m}\mathcal{J}_{\varepsilon}u^{\varepsilon})\right).
\end{multline}
We can then write $$II=II_{1}+II_{2}+II_{3}+II_{4},$$ where each of these terms corresponds to one of the
four terms on the right-hand side of (\ref{withProductRule}).  We will now deal with these one at a time.

We again will frequently use the fact that $\mathcal{J}_{\varepsilon}$ is self-adjoint and commutes
with $\partial_{x}.$  To begin, we have
$$II_{1}=\int_{0}^{2\pi}(\partial_{x}^{6}\mathcal{J}_{\varepsilon}u^{\varepsilon})
(\partial_{x}^{11}\mathcal{J}_{\varepsilon}u^{\varepsilon})\ dx.$$
After integrating by parts and recognizing a perfect derivative, we see that $II_{1}=0.$

For $II_{2},$ we have
$$II_{2}=-2\int_{0}^{2\pi}(\mathcal{J}_{\varepsilon}u^{\varepsilon})
(\partial_{x}^{6}\mathcal{J}_{\varepsilon}u^{\varepsilon})
(\partial_{x}^{7}\mathcal{J}_{\varepsilon}u^{\varepsilon})\ dx.$$
We recognize a perfect derivative and integrate by parts, finding
$$II_{2}=\int_{0}^{2\pi}(\mathcal{J}_{\varepsilon}u^{\varepsilon}_{x})
(\partial_{x}^{6}\mathcal{J}_{\varepsilon}u^{\varepsilon})^{2}\ dx\leq c\mathcal{E}^{3/2}.$$

For $II_{3},$ we have
$$II_{3}=\omega\int_{0}^{2\pi}(\partial_{x}^{6}\mathcal{J}_{\varepsilon}u^{\varepsilon})
(\partial_{x}^{7}\mathcal{J}_{\varepsilon}u^{\varepsilon})\ dx=0.$$
As before, this integrates to zero because it is the integral of a perfect derivative over a periodic
interval.

Finally, we treat $II_{4}.$  We have
\begin{equation}\label{leftoversToEstimate}
II_{4}=-2\sum_{m=1}^{6}{6\choose m}\int_{0}^{2\pi}
(\partial_{x}^{6}\mathcal{J}_{\varepsilon}u^{\varepsilon})
(\partial_{x}^{m}\mathcal{J}_{\varepsilon}u^{\varepsilon})
(\partial_{x}^{7-m}\mathcal{J}_{\varepsilon}u^{\varepsilon})\ dx.\end{equation}
Each of the integrals on the right-hand side of (\ref{leftoversToEstimate}) can be bounded by 
$\mathcal{E}^{3/2}.$ 

We therefore conclude that there exists $c>0$ such that
$$\frac{d\mathcal{E}}{dt}\leq c\mathcal{E}^{3/2}.$$
This can be rephrased as $$\frac{d\|u^{\varepsilon}\|_{H^{6}}}{dt}\leq c\|u^{\varepsilon}\|_{H^{6}}^{2}.$$
This estimate clearly indicates that the solutions $u^{\varepsilon}$ cannot blow up arbitrarily quickly,
and thus exist on a common time interval.  So, we have shown that there exists $T>0,$ independent
of $\varepsilon,$ such that for all $\varepsilon>0,$ we have 
$u^{\varepsilon}\in C([0,T]; H^{6}),$ with the norm bounded independently of $\varepsilon.$

Since $u^{\varepsilon}$ is bounded uniformly (in both $t$ and $\varepsilon$) and since
$u^{\varepsilon}$ solves (\ref{mollEvol}), we see that $\partial_{t}u^{\varepsilon}$ is uniformly
bounded in $H^{1},$ and thus in $L^{\infty}.$  
By the Arzela-Ascoli theorem, there exists a subsequence (which we do not relabel) 
and a limit, $u\in C([0,2\pi]\times[0,T]),$ such that $u_{\varepsilon}\rightarrow u$ in this space.
Standard arguments (again, cf. Chapter 3 of \cite{majdaBertozzi}, for instance) then imply
that $u$ belongs to the space $C([0,T]; H^{6}),$ that $u$ obeys the same uniform bound
as the $u^{\varepsilon},$ and that $u$ is a solution of the equation (\ref{fifthOrderEquation}).

Uniqueness of solutions (and, in fact, continuous dependence on the initial data) follows
from a more elementary version of the energy estimate.  If we let $u\in C([0,T];H^{6})$ be a solution 
corresponding to initial data $u_{0}\in H^{6},$ and if we let $v\in C([0,T];H^{6})$ be a solution 
corresponding to initial data $v_{0}\in H^{6},$ then we can estimate the norm of $u-v.$  If
$E_{d}=\|u-v\|_{L^{2}}^{2},$ then a straightforward calculation, together with the uniform bounds 
established previously, indicates $\frac{dE_{d}}{dt}\leq cE_{d}.$  This implies that $E_{d}\leq E_{d}(0)
e^{ct},$ for all $t$ for which the solutions are defined.  If $u_{0}=v_{0},$ then we see that $u=v.$
This is the desired uniqueness result.\hfill$\blacksquare$

\begin{remark} \label{uniformInOmega}
These estimates are uniform in $\omega,$ since all of the terms in the energy estimate
which involve $\omega$ are equal to zero.
\end{remark}

\subsection{Reformulation}\label{section:reformulation}

In this section, we use a normal form \cite{shatah} to rewrite the evolution equation in a beneficial way.

Taking the Fourier transform, we let $u_{k}$ be the 
Fourier coefficients of $u.$
We get the evolution equations
$$\partial_{t}u_{k}=ik^{5}u_{k}+i\omega ku_{k}-ik\sum_{j=-\infty}^{\infty}\primed u_{k-j}u_{j}.$$
We only consider $k\neq 0$ since we know already that $\partial_{t}u_{0}=u_{0}=0.$
The prime on the sum indicates that $j=0$ and $j=k$ are excluded, as these modes are unnecessary
since the mean of $u$ is equal to zero.
The initial condition is $$u_{k}(0)=\hat{g}(k).$$

We bring the linear terms to the left-hand side:
$$\partial_{t}u_{k}-ik^{5}u_{k}-i\omega ku_{k}=-ik\sum_{j=-\infty}^{\infty}\primed u_{k-j}u_{j}.$$
We use an integrating factor, so we define $v_{k}$ through the equation
$$v_{k}(t)=u_{k}(t)e^{-ik^{5}t-i\omega kt}.$$  This yields the following:
$$\partial_{t}v_{k}=-ik\sum_{j=-\infty}^{\infty}\primed
e^{-ik^{5}t-i\omega kt}e^{i(k-j)^{5}t+i\omega(k-j)t}e^{ij^{5}t+i\omega jt}v_{k-j}v_{j}.$$
The exponents simplify as the terms in the exponent related to the transport speed $\omega$ 
all cancel.  This leaves us with 
$$\partial_{t}v_{k}=-ik\sum_{j=-\infty}^{\infty}\primed
e^{-ik^{5}t}e^{i(k-j)^{5}t}e^{ij^{5}t}v_{k-j}v_{j}.$$
It will therefore be helpful to understand the quantity $k^{5}-(k-j)^{5}-j^{5},$ as this appears
in the exponent.

Using Pascal's triangle, we write
$$k^{5}=(k-j)^{5}+5(k-j)^{4}j+10(k-j)^{3}j^{2}+10(k-j)^{2}j^{3}+5(k-j)j^{4}+j^{5}.$$
Subtracting and factoring, we get
$$k^{5}-(k-j)^{5}-j^{5}=5(k-j)j((k-j)^{3}+2(k-j)^{2}j+2(k-j)j^{2}+j^{3}).$$
We add and subtract to make the final quantity more like $k^{3}:$
$$k^{5}-(k-j)^{5}-j^{5}=5(k-j)j(k^{3}-(k-j)^{2}j-(k-j)j^{2}).$$
This simplifies further:
$$k^{5}-(k-j)^{5}-j^{5}=5(k-j)jk(k^{2}-(k-j)j).$$
If we introduce the notation $\sigma=k^{2}-jk+j^{2},$ then this is
$$k^{5}-(k-j)^{5}-j^{5}=5(k-j)jk\sigma.$$
Using this identity, we are able to write the evolution equation for the $v_{k}:$
\begin{equation}\label{vkFirst}
\partial_{t}v_{k}=-ik\sum_{j=-\infty}^{\infty}\primed e^{-5i(k-j)jk\sigma t}v_{k-j}v_{j}.
\end{equation}

We can then rewrite this, recognizing that the exponential is in fact the derivative
of an exponential:
$$\partial_{t}v_{k}=-ik\sum_{j=-\infty}^{\infty}\primed
\left(\partial_{t}\left(\frac{e^{-5i(k-j)jk\sigma t}}
{-5i(k-j)jk\sigma}\right)\right)v_{k-j}v_{j}=\sum_{j=-\infty}^{\infty}\primed
\left(\partial_{t}\left(\frac{e^{-5i(k-j)jk\sigma t}}{5(k-j)j\sigma}\right)\right)v_{k-j}v_{j}.$$
Next, we ``differentiate by parts,'' moving the time derivative:
$$\partial_{t}v_{k}=\partial_{t}\left(\frac{1}{5}\sum_{j=-\infty}^{\infty}\primed
\frac{e^{-5i(k-j)jk\sigma t}}
{(k-j)j\sigma}
v_{k-j}v_{j}\right)-\frac{1}{5}\sum_{j=-\infty}^{\infty}\primed\frac{e^{-5i(k-j)jk\sigma t}}
{(k-j)j\sigma}\partial_{t}\left(v_{k-j}v_{j}\right),$$

We define $B$ through its Fourier coefficients, $B_{k}(t),$ as 
\begin{equation}\label{definitionOfB}
B_{k}(t)=-\frac{1}{5}\sum_{j=-\infty}^{\infty}\primed
\frac{e^{-5i(k-j)jk\sigma t}}{(k-j)j\sigma} v_{k-j}v_{j}.\end{equation}
We then are able to write the evolution equations as
$$\partial_{t}\left[v_{k}+B_{k}\right]=-\frac{1}{5}\sum_{j=-\infty}^{\infty}\primed\left[
\frac{e^{-5i(k-j)jk\sigma t}}
{(k-j)j\sigma}\right]\left[(\partial_{t}v_{k-j})v_{j}
+v_{k-j}(\partial_{t}v_{j})\right].$$
Next, we substitute for $\partial_{t}v_{k-j}$ and $\partial_{t}v_{j}.$
We let $\tilde{\sigma}=j^{2}-j\ell+\ell^{2};$ using (\ref{vkFirst}), we have the
following:
$$\partial_{t}v_{j}=-ij\sum_{\ell=-\infty}^{\infty}\primed
e^{-5itj\ell(j-\ell)\tilde{\sigma}}v_{j-\ell}v_{\ell}.$$
We are then able, using a symmetry between $k-j$ and $j,$ to write
\begin{equation}\label{doubleSum}
\partial_{t}[v_{k}+B_{k}]=\frac{2i}{5}\sum_{j=-\infty}^{\infty}\primed\sum_{\ell=-\infty}^{\infty}\primed
\left[\frac{e^{-5it[(k-j)jk\sigma+(j-\ell)j\ell\tilde{\sigma}]}}{(k-j)\sigma}\right]
v_{k-j}v_{j-\ell}v_{\ell}.\end{equation}  (We reiterate that $k\neq j.$)
We give the name $R_{k}$ to the right-hand side of
(\ref{doubleSum}), and we let $R$ be the function with Fourier coefficients equal to $R_{k}$ for
all $k.$
So, we have $\partial_{t}[v_{k}+B_{k}]=R_{k}.$  Integrating with respect to time, we have
\begin{equation}\label{integratedInTime}
v_{k}(t)-v_{k}(0)=B_{k}(0)-B_{k}(t)+\int_{0}^{t}R_{k}(s)\ ds.\end{equation}
We transform back to $u$ by multiplying (\ref{integratedInTime}) by $e^{ik^{5}t+i\omega kt},$ 
and we note
that $v_{k}(0)=u_{k}(0).$  These considerations yield the following:
\begin{equation}\label{ukt}u_{k}(t)-e^{ik^{5}t+i\omega kt}u_{k}(0)=e^{ik^{5}t+i\omega kt}\left(
B_{k}(0)-B_{k}(t)+\int_{0}^{t}R_{k}(s)\ ds\right).\end{equation}
Notice that (\ref{ukt}) is the $k^{\mathrm{th}}$ Fourier coefficient of the Duhamel integral at
time $t,$ or $\mathcal{F}(S_{D}(t)u_{0})(k).$

\subsection{Estimates}\label{section:estimates}

We now estimate $B$ and $R$ and associated quantities, to demonstrate the smoothing 
described in Theorem \ref{generalTheorem} for our equation with fifth-order dispersion.

\begin{remark} \label{uniformInOmega2}
It will be plain to see that all estimates made in the present section
are uniform in $\omega.$
\end{remark}

\begin{lemma}\label{BLemma}
If $s\geq 1$ and $v\in H^{s},$ then $B\in H^{s+3},$ with the estimate
$$\|B\|_{H^{s+3}}^{2}\leq c \|v\|_{H^{s}}^{4}.$$
\end{lemma}
(We note that for our main theorem, 
we only actually need $B\in H^{s+2},$ but it turns out that $B\in H^{s+3}.$)

\noindent{\bf Proof:} We will show that $\partial_{x}^{3}B$ is in $H^{s}.$
We begin by taking three derivatives of $B,$ which requires multiplying
(\ref{definitionOfB}) by $(ik)^{3}:$
$$(ik)^{3}B_{k}=\frac{i}{5}\sum_{j=-\infty}^{\infty}\primed E(k,j)
\frac{k^{3}}{j(k-j)(k^{2}-kj+j^{2})}v_{k-j}v_{j},$$
where $E(k,j)$ represents the exponential,
$E(k,j)=e^{-5itjk(k-j)\sigma}.$

We will demonstrate now that $\displaystyle\frac{k^{3}}{j(k-j)(k^{2}-kj+j^{2})}$ is bounded by a constant.
To begin, we consider $$\left|\frac{k}{j(k-j)}\right|=\left|\frac{k-j+j}{j(k-j)}\right|\leq
\left|\frac{1}{j}\right|+\left|\frac{1}{k-j}\right|\leq 2.$$  Next, we consider $\sigma=k^{2}-kj+j^{2}.$
We observe that 
\begin{equation}\label{sigmaForm}
\sigma=\frac{1}{2}k^{2}+\frac{1}{2}(k-j)^{2}+\frac{1}{2}j^{2},
\end{equation}
so clearly $\sigma\geq\frac{1}{2}k^{2}.$ Thus,
\begin{equation}\label{sigmaGain}
\left|\frac{k^{2}}{\sigma}\right|
\leq 2.
\end{equation}
Thus, for any $k$ and any $j,$ we have
\begin{equation}\label{boundBy4}
\left|E(k,j)\frac{k^{3}}{j(k-j)(k^{2}-kj+j^{2})}\right|\leq 4.
\end{equation}
We give the name $\Phi(k,j)=E(k,j)\frac{k^{3}}{j(k-j)(k^{2}-kj+j^{2})}.$

Of course, we have $$\|\partial_{x}^{s+3}B\|_{L^{2}}^{2}=\||k|^{s}(ik^{3}B_{k})\|_{\ell^{2}}^{2}.$$
We then have
$$\||k|^{s}(ik^{3}B_{k})\|_{\ell^{2}}^{2}=\frac{1}{25}
\sum_{k=-\infty}^{\infty}\sum_{j=-\infty}^{\infty}\primed
\sum_{\ell=-\infty}^{\infty}\primed
\Phi(k,j)\bar{\Phi}(k,\ell)k^{2s}v_{k-j}v_{j}\bar{v}_{k-\ell}\bar{v}_{\ell}.$$
In light of (\ref{boundBy4}), we have
\begin{equation}\label{convolutionsAndAbsolutes}
\||k|^{s}(ik^{3}B_{k})\|_{\ell^{2}}^{2}\leq 
\sum_{k,j,\ell}k^{2s}|v_{k-j}||v_{j}||v_{k-\ell}||v_{\ell}|.\end{equation}

Let $V$ be the function defined through its Fourier coefficients as $\mathcal{F}V(k)=|v_{k}|,$ 
for all $k.$  Note that since $v\in H^{s},$ we have $V\in H^{s},$ with $\|v\|_{H^{s}}=\|V\|_{H^{s}}.$
Since $V\in H^{s}$ with $s\geq 1,$ 
we can see that $V^{2}\in H^{s},$ with $\|V^{2}\|_{H^{s}}\leq c\|v\|_{H^{s}}^{2}.$
Notice that the right-hand side of (\ref{convolutionsAndAbsolutes}) is equal to 
$\|\partial_{x}^{s}(V^{2})\|_{L^{2}}^{2}.$
This completes the proof. \hfill $\blacksquare$

The particular estimate we need for $B$ follows from Lemma \ref{BLemma}:

\begin{cor}\label{bCorollary}
 If $s\geq 1$ and $u\in C([0,T]; H^{s}),$ and if $t\in[0,T],$ then
$$\mathcal{F}^{-1}\left(e^{ik^{5}t+i\omega kt}(B_{k}(0)-B_{k}(t))\right)\in H^{s+3},$$ with the bound
$$\left\|\mathcal{F}^{-1}\left(e^{ik^{5}t+i\omega kt}(B_{k}(0)-B_{k}(t))\right)\right\|_{H^{s+3}}\leq 
c\|u\|_{C([0,T]; H^{s})}^{2}.$$
\end{cor}

\noindent{\bf Proof:} This follows immediately from Lemma \ref{BLemma}, and from the fact that
$|e^{i\theta}|=1$ for any real $\theta.$
\hfill$\blacksquare$

Having proved a satisfactory estimate for $B,$ we turn to $R.$

\begin{lemma}\label{RLemma}
If $s\geq 1$ and $v\in H^{s},$ then $R\in H^{s+2},$ with the estimate
$$\|R\|_{H^{s+2}}^{2}\leq c \|v\|_{H^{s}}^{6}.$$
\end{lemma}

\noindent{\bf Proof:}  Recall the formula for $R_{k},$
$$R_{k}=\frac{2i}{5}\sum_{j=-\infty}^{\infty}\!\!\!{}^{'}
\sum_{\ell=-\infty}^{\infty}\!\!\!{}^{'}\frac{1}{(k-j)\sigma}
(v_{k-j}v_{j-\ell}v_{\ell})\exp\{(-5it)(kj(k-j)\sigma+j\ell(j-\ell)\tilde{\sigma})\},
$$
with $\sigma=k^{2}-kj+j^{2}$ and $\tilde{\sigma}=j^{2}-j\ell+\ell^{2}.$
As we noted in the proof of Lemma \ref{BLemma}, we have $\frac{k^{2}}{\sigma}\leq 2.$
Following the lines of the proof of Lemma \ref{BLemma}, we arrive at
$$\|\partial_{x}^{s+2}R\|_{L^{2}}^{2}\leq c\sum_{k,j,\ell,m,n}k^{2s}
|v_{k-j}||v_{j-\ell}||v_{\ell}||v_{k-m}||v_{m-n}||v_{n}|.$$
Letting $V$ be as in the proof of Lemma \ref{BLemma}, we see that the right-hand side
is a constant times square of  the $L^{2}$ norm of $\partial_{x}^{s}(V^{3}).$  Thus, this is bounded by
$\|v\|_{H^{s}}^{6},$ as claimed. \hfill$\blacksquare$

Lemma \ref{RLemma} implies the following, which is the estimate we need for $R:$

\begin{cor}\label{rCorollary}
If $s\geq 1$ and $u\in C([0,T];H^{s}),$ and if $t\in[0,T],$ then
$$\mathcal{F}^{-1}\left(e^{ik^{5}t+i\omega kt}\int_{0}^{t}R_{k}(s)\ ds\right)\in H^{s+2},$$ with
$$\left\|\mathcal{F}^{-1}\left(e^{ik^{5}t+i\omega kt}\int_{0}^{t}R_{k}(s)\ ds\right)\right\|_{H^{s+2}}\leq 
c\|u\|_{C([0,T]; H^{s})}^{3}.$$
\end{cor}

\noindent{\bf Proof:} We begin by noting that, of course,
\begin{multline}\nonumber
\left\|\partial_{x}^{s+2}\mathcal{F}^{-1}\left(e^{ik^{5}t+i\omega kt}
\int_{0}^{t}R_{k}(s)\ ds\right)\right\|_{L^{2}}^{2}
=\left\|(ik)^{s+2}e^{ik^{5}t+i\omega kt}\int_{0}^{t}R_{k}(s)\ ds\right\|_{\ell^{2}}^{2}\\
=\left\|\int_{0}^{t}k^{s+2}R_{k}(s)\ ds\right\|_{\ell^{2}}^{2}
=\sum_{k=-\infty}^{\infty}\left(\int_{0}^{t}k^{s+2}R_{k}(s)\ ds\right)
\left(\int_{0}^{t}k^{s+2}\bar{R}_{k}(\tau)\ d\tau\right).
\end{multline}
We use the triangle inequality:
\begin{multline}\nonumber
\left\|\partial_{x}^{s+2}\mathcal{F}^{-1}\left(e^{ik^{5}t+i\omega kt}
\int_{0}^{t}R_{k}(s)\ ds\right)\right\|_{L^{2}}^{2}\\
\leq\sum_{k=-\infty}^{\infty}\left(\int_{0}^{t}|k^{s+2}R_{k}(s)|\ ds\right)
\left(\int_{0}^{t}|k^{s+2}R_{k}(\tau)|\ d\tau\right).
\end{multline}
By Tonelli's Theorem, we can exchange the sum and the integrals:
\begin{multline}\nonumber
\left\|\partial_{x}^{s+2}\mathcal{F}^{-1}\left(e^{ik^{5}t+i\omega kt}
\int_{0}^{t}R_{k}(s)\ ds\right)\right\|_{L^{2}}^{2}\\
\leq \int_{0}^{t}\int_{0}^{t}\sum_{k=-\infty}^{\infty}|k^{s+2}R_{k}(s)||k^{s+2}R_{k}(\tau)|\ ds\ d\tau.
\end{multline}
We then use the Cauchy-Schwartz inequality:
\begin{multline}\nonumber
\left\|\partial_{x}^{s+2}\mathcal{F}^{-1}\left(e^{ik^{5}t+i\omega kt}\int_{0}^{t}R_{k}(s)\ ds\right)\right\|_{L^{2}}^{2}\\
\leq c\int_{0}^{t}\int_{0}^{t}\|R(\cdot,s)\|_{H^{s+2}}\|R(\cdot,\tau)\|_{H^{s+2}}\ ds\ d\tau.
\end{multline}
We then use Lemma \ref{RLemma}, and the proof is complete.\hfill$\blacksquare$

We are now in a position to show that the necessary estimate for the Duhamel integral
holds for the equation (\ref{fifthOrderEquation}).  
\begin{theorem} \label{smoothingEstimate}
Let $0<T_{1}<T_{2}$ be given.  There exists $\gamma>0$ such that for any
$u_{0}\in H^{6}$ such that $\|u_{0}\|_{H^{6}}<\gamma,$ then there is a unique solution of the initial value
problem (\ref{fifthOrderEquation}) with initial data $u_{0},$ with the solution $u\in C([0,T_{2}];H^{6}).$
There exists $c>0$ and $\tilde{\gamma}\in(0,\gamma)$ such that
for any $T\in[T_{1},T_{2}],$ and for any $u_{0}\in H_{0}^{6}$ such that 
$\|u_{0}\|_{H^{6}}<\tilde{\gamma},$ then
$$\|S_{D}(T)u_{0}\|_{H^{8}}\leq c\|u_{0}\|_{H^{6}}^{2}.$$
\end{theorem}

\noindent{\bf{Proof:}}
The formula (\ref{ukt}) and the estimates of Corollary \ref{bCorollary} and Corollary \ref{rCorollary} 
for $B$ and $R,$ respectively, immediately imply
$$\|S_{D}(T)u_{0}\|_{H^{8}}\leq c\|u\|_{C([0,T];H^{6})}^{2}.$$
However, we are not yet finished because we need the bound to be in terms of the initial data, and not 
in terms of the solution at positive times.

As discussed in Section \ref{section:existenceTheory} above, we have
$$\frac{d}{dt}\|u^{\varepsilon}\|_{H^{6}}\leq c\|u^{\varepsilon}\|_{H^{6}}^{2}.$$
Let $\|u_{0}\|_{H^{6}}=\frac{\delta}{2}.$  Then, as long as $\|u^{\varepsilon}(\cdot,t)\|_{H_{6}}\leq\delta=
2\|u_{0}\|_{H^{6}},$
we have $$\frac{d}{dt}\|u^{\varepsilon}\|_{H^{6}}\leq c\delta \|u^{\varepsilon}\|_{H^{6}},$$
and thus $$\|u^{\varepsilon}\|_{H^{6}}\leq \|u_{0}\|_{H^{6}}e^{c\delta t}=\frac{\delta e^{c\delta t}}{2}.$$
This implies that $\|u^{\varepsilon}\|_{H^{6}}\leq 2\|u_{0}\|_{H^{6}}$ as long as $e^{c\delta t}<2.$  This is 
valid as long as $t<\ln(2)/c\delta;$ notice that this bound goes to infinity as $\delta$ vanishes (that is,
the ``doubling time'' for solutions goes to infinity as the initial size of the solutions goes to zero).
Taking the limit as $\varepsilon$ vanishes, then (along our subsequence), we find
$$\|u\|_{H^{6}}\leq 2 \|u_{0}\|_{H^{6}},$$ as long as $t<\ln(2)/c\delta.$

Given the set of potential temporal periods of interest,  $[T_{1},T_{2}],$ we may choose $\delta$
sufficiently small so that as long as $\|u_{0}\|_{H^{6}}<\frac{\delta}{2},$ then
for all $t\in[0,T_{2}],$ we have $\|u(\cdot,t)\|_{H^{6}}<2\|u_{0}\|_{H^{6}}.$
We have shown above that $\|S_{D}(T)u_{0}\|_{H^{8}}\leq c \|u\|_{C([0,T_{2}]; H^{6})}^{2}.$  We now
can bound this in terms of $\|u_{0}\|_{H^{6}},$ so that
$$\|S_{D}(T)u_{0}\|_{H^{8}}\leq 4c\|u_{0}\|_{H^{6}}^{2},$$
for any $T\in[T_{1},T_{2}],$ and for any $u_{0}$ satisfying our smallness assumption.  This is the 
desired bound.  \hfill$\blacksquare$

\subsection{Completion of the example}

In this section, we state a specific theorem on nonexistence of time-periodic solutions, making use
of the above.  Consider the equation
\begin{equation}\label{specificEquation}
\partial_{t}\tilde{u}=\partial_{x}^{5}\tilde{u}-2\tilde{u}\partial_{x}\tilde{u}.
\end{equation}
Consider the initial data, $\tilde{u}_{0}\in H^{6},$ with the mean of $\tilde{u}_{0}$ equal to $\alpha.$
As above, we define $u_{0}=\tilde{u}_{0}-\alpha,$ and we let $u=\tilde{u}-\alpha.$
The evolution equation satisfied by $u,$ as discussed above, is
\begin{equation}\label{specificModified}
\partial_{t}u=\partial_{x}^{5}u-2\alpha\partial_{x}u-2u\partial_{x}u.
\end{equation}
If we let $A=\partial_{x}^{5}-2\alpha\partial_{x},$ then we see that the symbol of $A$ is
\begin{equation}\label{specificSymbol}
\hat{A}(k)=i(k^{5}-2\alpha k)=ik(k^{4}-2\alpha).
\end{equation}
Thus, if $2\alpha<1,$ then there are no zeros of the symbol in $\mathbb{Z}\setminus\{0\}.$
Recalling the hypotheses {\bf(H)}, we have $M=2,$ and we 
let $r_{1}=5,$ $r_{2}=1,$ $Z_{1}=\{1\},$ and
$Z_{2}=(-\frac{1}{2},\frac{1}{2}).$  
Letting $\alpha_{1}\in Z_{1}$ and $\alpha_{2}\in Z_{2},$ we see that we may take
$\beta_{1}=\frac{1}{2}$ and $\beta_{2}=\frac{1}{2}.$  Letting $\alpha_{2}=-2\alpha,$ we see that
the hypotheses {\bf(H)} are satisfied, and Lemma \ref{linearLemma2} and Corollary \ref{almostEveryCorollary}
hold, with uniform estimates for $\alpha\in(-\frac{1}{4},\frac{1}{4}).$
Therefore, Theorem \ref{generalTheorem} applies for $\alpha\in(-\frac{1}{4},\frac{1}{4}).$
Theorem \ref{smoothingEstimate} also applies, and in light of Remarks \ref{uniformInOmega}
and \ref{uniformInOmega2}, we see that the constants in Theorem \ref{smoothingEstimate} can be taken
to be uniform with respect to $\alpha.$
This implies that equation (\ref{specificModified}) does not possess small, nonzero time-periodic
solutions, uniformly in $\alpha\in(-\frac{1}{4},\frac{1}{4}).$  

Adding the mean $\alpha$ back to $u,$
we get $\tilde{u}=(\tilde{u}-\alpha)+\alpha,$ with $u=\tilde{u}-\alpha.$ We know that for 
$\alpha\in(-\frac{1}{4},\frac{1}{4}),$ $\tilde{u}-\alpha$
does not possess small, nonzero time-periodic solutions with the associated periods ($T\in W_{p,\delta}$ or
$T\in W_{p},$ as appropriate).  Furthermore, we know that $\|\tilde{u}\|_{H^{6}}\geq \|\tilde{u}-\alpha\|_{H^{6}}.$
This implies that the only small time-periodic solutions of (\ref{specificEquation}) with the
given temporal periods are $\tilde{u}=\alpha.$  Thus
(\ref{specificEquation}) does not possess small, non-constant time-periodic solutions with 
the given temporal periods.  This proves the following corollary:

\begin{cor}\label{mainSpecificResult}
Let $0<T_{1}<T_{2}$ be given.  Let $p\in(1,2]$ be given.  Let $0<\delta<T_{2}-T_{1}$ be
given.  
Let $W_{p,\delta}\subseteq[T_{1},T_{2}]$ be as in Lemma \ref{linearLemma2},
with $A$ given by (\ref{specificSymbol}), for $|\alpha|\leq\frac{1}{4}.$
There exists $r_{1}>0$ such that 
for all $T\in W_{p,\delta},$ if
$u$ is a smooth,
non-constant time-periodic solution of (\ref{specificEquation}) with temporal period $T,$ then 
$$\inf_{t\in[0,T]}\|u\|_{H^{6}}>r_{1}.$$

Let $W_{p}\subseteq[T_{1},T_{2}]$ be as in Corollary \ref{almostEveryCorollary}.  
Let $T\in W_{p}$ be given.  Then there exists
$r_{2}$ such that if $u$ is a smooth, non-constant time-periodic solution of (\ref{specificEquation}) 
with temporal period $T,$ then
$$\inf_{t\in[0,T]}\|u\|_{H^{6}}>r_{2}.$$
\end{cor}

\section{Further examples with $\tilde{p}=0$}\label{section:further}

In this section, we provide a few other equations which can be treated similarly to the above.
We do not provide full proofs in this section, but instead point out the differences with the prior 
proof.

\subsection{Non-resonant Kawahara equations}\label{section:kawahara}

The Kawahara equation has been justified as a model for water waves with surface tension
\cite{duell},  \cite{schneiderWayne-ARMA}.  It can be written as 
\begin{equation}\label{specificKawahara}
\partial_{t}\tilde{u}=\partial_{x}^{5}\tilde{u}-\theta\partial_{x}^{3}\tilde{u}-2\tilde{u}\partial_{x}\tilde{u},
\end{equation}
with $\theta>0.$  As before, we take this with initial data $\tilde{u}(\cdot,0)=\tilde{u}_{0},$ and we
assume that the mean of $\tilde{u}_{0}$ is equal to $\alpha.$  We again let $u=\tilde{u}-\alpha,$
and we find that the equation satisfied by $u$ is
\begin{equation}\label{zeroMeanKawahara}
\partial_{t}u=\partial_{x}^{5}u-\theta\partial_{x}^{3}u-2\alpha\partial_{x}u-2u\partial_{x}u.
\end{equation}

Our prior results extend to the Kawahara equation as long as the constant $\theta$ is chosen
to avoid resonance.  In particular, we must require $k^{5}-\theta k^{3}\neq0,$ for all 
$k\in\mathbb{Z}\setminus\{0\}.$  Notice that this is the same as requiring
\begin{equation}\label{nonresonantCondition}
\min_{k\in\mathbb{Z}\setminus\{0\}}\left|k^{5}-\theta k^{3}\right|>0.
\end{equation}
This implies that there exists a constant $\bar{\alpha}>0$ and a constant $\beta_{2}>0$
such that for all $\alpha\in(-\bar{\alpha},\bar{\alpha}),$ for all $k\in\mathbb{Z}\setminus\{0\},$
$$\left|k^{5}-\theta k^{3}-2\alpha k\right|\geq \beta_{2}.$$  
We now verify {\bf(H)}, taking $M=3.$  We let $r_{1}=5,$ $r_{2}=3,$ and $r_{3}=1.$
We let $Z_{1}=\{1\}$ and $Z_{2}=\{\theta\}.$  We take $\alpha_{1}=1,$ $\alpha_{2}=\theta,$
and $\alpha_{3}=-2\alpha,$ with $Z_{3}=(-2\bar{\alpha},2\bar{\alpha}).$
Then, {\bf(H)} is satisfied, with $\beta_{1}=\frac{1}{2}.$

This means that Theorem \ref{generalTheorem} applies, and all that must be done to conclude
the nonexistence of small doubly periodic waves for the nonresonant 
Kawahara equation is that the smoothing property must be demonstrated.
We are able to demonstrate the smoothing property with $\tilde{p}=0.$

We take the Fourier transform of (\ref{zeroMeanKawahara}), finding
$$\partial_{t}u_{k}=ik^{5}u_{k}+i\theta k^{3}u_{k}-2i\alpha ku_{k}
-ik\sum_{j=-\infty}^{\infty}u_{k-j}u_{j}.$$
As before, we use an integrating factor, defining 
$$v_{k}=u_{k}\exp\{it(-k^{5}-\theta k^{3}+2\alpha k)\}.$$
We have the following evolution equation for $v_{k}:$
\begin{equation}\nonumber
\partial_{t}v_{k}=-ik\sum_{j=-\infty}^{\infty}\primed(v_{k-j}v_{j})
\exp\{(it)\Phi(j,k)\},
\end{equation}
where the phase function, $\Phi,$ is given by
\begin{equation}\nonumber
\Phi(j,k)=-k^{5}-\theta k^{3}+2\alpha k+(k-j)^{5}+\theta(k-j)^{3}-2\alpha(k-j)
+j^{5}+\theta j^{3}-2\alpha j.
\end{equation}
This simplifies, as all the terms with $\alpha$ cancel, and also since we have previously
computed $k^{5}-(k-j)^{5}-j^{5}.$  We have not previously computed $k^{3}-(k-j)^{3}-j^{3},$
but this is straightforward: $$k^{3}-(k-j)^{3}-j^{3}=3jk(k-j).$$
These considerations imply the following:
\begin{equation}\nonumber
\Phi(j,k)=5jk(k-j)\sigma+3\theta jk(k-j)=5jk(k-j)\left(\sigma+\frac{3}{5}\theta\right).
\end{equation}

The critical step in the proof of smoothing in Section \ref{section:fifthOrder} was
inequality (\ref{sigmaGain}).  We see that the corresponding inequality in the
present case is 
$$\left|\frac{k^{2}}{\sigma+\frac{3}{5}\theta}\right|\leq 2,$$
which follows immediately from (\ref{sigmaForm}) and the condition $\theta>0.$  
The rest of the proof of Section \ref{section:fifthOrder}
can be repeated, establishing the following:
\begin{cor}\label{kawaharaSpecificResult}
Let $0<T_{1}<T_{2}$ be given.  Let $p\in(1,2]$ be given.  Let $0<\delta<T_{2}-T_{1}$ be
given.  Let $\theta>0$ satisfy (\ref{nonresonantCondition}), and let $\bar{\alpha}>0$ be as above.
Let $W_{p,\delta}\subseteq[T_{1},T_{2}]$ be as in Lemma \ref{linearLemma2},
with $A$ given by $\mathcal{F}(A)=i(-k^{5}-\theta k^{3}+2\alpha k),$ for $|\alpha|\leq\bar{\alpha}.$
There exists $r_{1}>0$ such that 
for all $T\in W_{p,\delta},$ if
$u$ is a smooth, non-constant time-periodic solution of (\ref{specificKawahara}) 
with temporal period $T,$ then $$\inf_{t\in[0,T]}\|u\|_{H^{6}}>r_{1}.$$

Let $W_{p}\subseteq[T_{1},T_{2}]$ be as in Corollary \ref{almostEveryCorollary}.  
Let $T\in W_{p}$ be given.  Then there exists
$r_{2}$ such that if $u$ is a smooth, non-constant time-periodic solution of 
(\ref{specificKawahara}) 
with temporal period $T,$ then
$$\inf_{t\in[0,T]}\|u\|_{H^{6}}>r_{2}.$$
\end{cor}

\begin{remark} Above, we appeared to use in a fundamental way the property $\theta>0.$
We assume $\theta>0$ only because this appears to be a feature of the Kawahara equation
as it exists in the prior literature.  If we instead had $\theta<0,$ our argument would still work.
For a particular value of $\theta,$ if there exist $(j,k)\in (\mathbb{Z}\setminus\{0\})^{2}$ such that
$\sigma+\frac{3}{5}\theta=0,$ then we treat such values of $(j,k)$ differently.
The arguments of Section \ref{section:fifthOrder} continue to apply whenever this quantity does
not vanish.  All that remains is to observe that, as can be seen from (\ref{sigmaForm}),
the set of $(j,k)$ for which $\sigma+\frac{3}{5}\theta$ does vanish is bounded for a fixed
value of $\theta$ (or for values of $\theta$ in a bounded set), and that regularity is determined by behavior for large $k.$
\end{remark}

\subsection{Seventh-order equations}\label{section:seventhOrder}

The proof of smoothing for the fifth-order equation above depended on combinatorial 
properties,
and the use of Pascal's triangle in particular.  With dispersion of seventh order, 
we must calculate $k^{7}-(k-j)^{7}-j^{7}.$

To begin, Pascal's triangle yields the following:
\begin{multline}\nonumber
k^{7}-(k-j)^{7}-j^{7}
=7\Big((k-j)^{6}j\\
+3(k-j)^{5}j^{2}+5(k-j)^{4}j^{3}+5(k-j)^{3}j^{4}+3(k-j)^{2}j^{5}+(k-j)j^{6}\Big).
\end{multline}
We can factor out $(k-j)j:$
\begin{multline}\nonumber
k^{7}-(k-j)^{7}-j^{7}
=7(k-j)j\Big((k-j)^{5}\\
+3(k-j)^{4}j+5(k-j)^{3}j^{2}+5(k-j)^{2}j^{3}+3(k-j)j^{4}+j^{5}\Big).
\end{multline}
We add and subtract, to make the quantity in parentheses more like $k^{5}:$
\begin{equation}\nonumber
k^{7}-(k-j)^{7}-j^{7}
=7(k-j)j\Big(k^{5}
-2(k-j)^{4}j-5(k-j)^{3}j^{2}-5(k-j)^{2}j^{3}-2(k-j)j^{4}\Big).
\end{equation}
From most of the terms in parentheses, then, we can again factor out $(k-j)j:$
\begin{equation}\nonumber
k^{7}-(k-j)^{7}-j^{7}
=7(k-j)j\Big(k^{5}
-(k-j)j\Big(2(k-j)^{3}+5(k-j)^{2}j+5(k-j)j^{2}+2j^{3}\Big)\Big).
\end{equation}
We again add and subtract, this time to produce $2k^{3}:$
\begin{equation}\nonumber
k^{7}-(k-j)^{7}-j^{7}
=7(k-j)j\Big(k^{5}-(k-j)j\Big(2k^{3}-k(k-j)j\Big)\Big).
\end{equation}
We can then factor out $k:$
\begin{equation}
k^{7}-(k-j)^{7}-j^{7}=7(k-j)jk\tau,
\end{equation}
where $$\tau=k^{4}-2k^{2}(k-j)j+(k-j)^{2}j^{2}.$$

We multiply $\tau$ out:
$$\tau=k^{4}-2k^{3}j+3k^{2}j^{2}-2kj^{3}+j^{4}.$$
We calculate the following:
$$\frac{1}{2}(k-j)^{4}
=\frac{1}{2}k^{4}-2k^{3}j+3k^{2}j^{2}-2kj^{3}+\frac{1}{2}j^{4}.$$
Therefore, we have
$$\tau=\frac{1}{2}k^{4}+\frac{1}{2}(k-j)^{4}+\frac{1}{2}j^{4};$$
note the similarity to (\ref{sigmaForm}).
Clearly, we have the inequality
\begin{equation}\label{gains4Derivatives}
\left|\frac{k^{4}}{\tau}\right|=\frac{k^{4}}{\tau}\leq 2.
\end{equation}

We can then perform all of the calculations of Section \ref{section:fifthOrder} for the equation
\begin{equation}\label{seventhOrderSpecificExample}
\partial_{t}u=\partial_{x}^{7}u-2u\partial_{x}u.
\end{equation}
Recall that $B$ and $R,$ defined in Section \ref{section:reformulation}, 
were shown in Section \ref{section:estimates} to gain two derivatives.  This property hinged on the
inequality (\ref{sigmaGain}).
In the present setting, the analogues of $B$ and $R$ would now gain four derivatives because of 
(\ref{gains4Derivatives}).  (As noted in the introduction, this allows one to see that the smoothing
mechanism we are using
is different than Kato smoothing, as Kato smoothing would provide a gain of three
derivatives with seventh-order dispersion.)  Following the argument of Section \ref{section:fifthOrder},
but using (\ref{gains4Derivatives}), we arrive at the following:

\begin{cor}\label{seventhOrderSpecificResult}
Let $0<T_{1}<T_{2}$ be given.  Let $p\in(1,4]$ be given.  Let $0<\delta<T_{2}-T_{1}$ be
given.  
Let $W_{p,\delta}\subseteq[T_{1},T_{2}]$ be as in Lemma \ref{linearLemma2},
with $A$ given by $\mathcal{F}(A)=i(-k^{7}+\alpha k),$ for $|\alpha|\leq\frac{1}{2}.$
There exists $r_{1}>0$ such that 
for all $T\in W_{p,\delta},$ if
$u$ is a smooth, non-constant time-periodic solution of (\ref{seventhOrderSpecificExample}) 
with temporal period $T,$ then $$\inf_{t\in[0,T]}\|u\|_{H^{8}}>r_{1}.$$

Let $W_{p}\subseteq[T_{1},T_{2}]$ be as in Corollary \ref{almostEveryCorollary}.  
Let $T\in W_{p}$ be given.  Then there exists
$r_{2}$ such that if $u$ is a smooth, non-constant time-periodic solution of 
(\ref{seventhOrderSpecificExample}) 
with temporal period $T,$ then
$$\inf_{t\in[0,T]}\|u\|_{H^{8}}>r_{2}.$$
\end{cor}

We note that a change from Corollary \ref{seventhOrderSpecificResult} as compared to
Corollary \ref{mainSpecificResult} is that we are now using $H^{8}$ instead of $H^{6}.$
The function space is chosen so that the solutions under consideration are classical solutions;
using the space $H^{8},$ the seventh derivative of $u$ appearing on the
right-hand side of (\ref{seventhOrderSpecificExample}) is classically defined.

\section{The KdV equation}\label{section:kdv}

We study the equation 
\begin{equation}\label{kdv1}
\partial_{t}\tilde{u}=-\partial_{x}^{3}\tilde{u}-\tilde{u}\partial_{x}\tilde{u}.
\end{equation}
(This is not the most traditional choice of coefficients for the KdV equation,
but the coefficients are changeable by scaling, so it makes no difference.)
The evolution equation (\ref{kdv1}) is taken with initial data 
$u(\cdot,0)=\tilde{g}.$
We let the mean of $\tilde{g}$ be denoted as $\alpha,$ and we define 
$g=\tilde{g}-\alpha.$  
Then, noticing that (\ref{kdv1})
preserves the mean of the solution, we define $u=\tilde{u}-\alpha,$ so that
the mean of $u$ is equal to zero, as long as the solution 
$\tilde{u}$ of (\ref{kdv1})
exists.  The evolution equation satisfied by $u$ is
\begin{equation}\label{kdv2}
\partial_{t}u=-\partial_{x}^{3}u-\alpha\partial_{x}u-u\partial_{x}u.
\end{equation}

We note that there are very many papers in the literature treating the well-posedness
of the KdV equation.  For example, in \cite{bourgainKdVGWP}, global well-posedness for
the periodic KdV equation with initial data in $H^{s}$ for $s\geq 0$ is established, and 
in  \cite{kappelerTopalov}, global well-posedness is established in $H^{s}$ for $s\geq -1.$
Nevertheless, we remark that the simple well-posedness proof for (\ref{fifthOrderEquationTilde}) given
in Section \ref{section:existenceTheory} 
can be suitably and straightforwardly modified to the KdV equation,
yielding the same results: that classical solutions exist for arbitrarily long intervals of time
if the initial data is sufficiently small, and that the doubling time of such solutions goes to
infinity as the size of the data vanishes.    

We take the Fourier transform of (\ref{kdv2}):
\begin{equation}\nonumber
\partial_{t}u_{k}=i(k^{3}-\alpha k)u_{k}
-i\sum_{j=-\infty}^{\infty}\primed u_{k-j}(ju_{j}).
\end{equation}
Since the mean of $u$ is equal to zero, the values $j=0$ and $j=k$ are 
unnecessary; as before, the prime indicates that these indices are excluded from the
summation.
We use an integrating factor, defining 
$v_{k}=u_{k}\exp\{-it(k^{3}-\alpha k)\}:$
\begin{equation}\nonumber
\partial_{t}v_{k}=-i\sum_{j=-\infty}^{\infty}\primed
\exp\{-it(k^{3}-(k-j)^{3}-j^{3})\}v_{k-j}(jv_{j}).
\end{equation}
(Notice that the terms in the exponential involving $\alpha$ all canceled.)
From Pascal's triangle, we see that $k^{3}-(k-j)^{3}-j^{3}=3kj(k-j).$  Thus,
we may write the following:
\begin{equation}\label{kdvVkBeforeDivision}
\partial_{t}v_{k}=-i\sum_{j=-\infty}^{\infty}\primed
\exp\{-3itkj(k-j)\}v_{k-j}(jv_{j}).
\end{equation}

We now manipulate the exponential:
\begin{equation}\nonumber
\partial_{t}v_{k}=-i\sum_{j=-\infty}^{\infty}\primed
\partial_{t}\left(\frac{e^{-3itkj(k-j)}}{-3ikj(k-j)}\right)
v_{k-j}(jv_{j}).
\end{equation}
We cancel the factor of $-i,$ and we also ``differentiate by parts:''
\begin{equation}\nonumber
\partial_{t}v_{k}=\partial_{t}\left(\sum_{j=-\infty}^{\infty}\primed
\frac{e^{-3itkj(k-j)}}{3kj(k-j)}v_{k-j}(jv_{j})\right)
-\sum_{j=-\infty}^{\infty}\frac{e^{-3itkj(k-j)}}{3kj(k-j)}\partial_{t}\left(
v_{k-j}(jv_{j})\right).
\end{equation}
We can write this as 
\begin{equation}\nonumber
\partial_{t}(v_{k}+B_{k})=R_{k},
\end{equation}
with 
\begin{equation}\nonumber
B_{k}=-\frac{1}{3}\sum_{j=-\infty}^{\infty}
\frac{e^{-3itkj(k-j)}}{kj(k-j)}v_{k-j}(jv_{j}),
\end{equation}
\begin{equation}\label{kdvRk1}
R_{k}=-\frac{1}{3}\sum_{j=-\infty}^{\infty}\primed
\frac{e^{-3itkj(k-j)}}{kj(k-j)}\partial_{t}\left(v_{k-j}(jv_{j})\right).
\end{equation}
We continue to rewrite this, by applying the time derivative on the right-hand side
of (\ref{kdvRk1}).
Using (\ref{kdvVkBeforeDivision}), we have
\begin{equation}\nonumber
\partial_{t}v_{k-j}=-i\sum_{m=-\infty}^{\infty}\primed
\exp\{-3it(k-j)m(k-j-m)\}v_{k-j-m}(mv_{m}),
\end{equation}
\begin{equation}\nonumber
\partial_{t}v_{j}=-i\sum_{\ell=-\infty}^{\infty}\primed
\exp\{-3itj\ell(j-\ell)\}v_{j-\ell}(\ell v_{\ell}).
\end{equation}
Using these, we write $R_{k}=R^{1}_{k}+R^{2}_{k},$ with
\begin{equation}\nonumber
R^{1}_{k}=\frac{i}{3}\sum_{j=-\infty}^{\infty}\primed
\sum_{m=-\infty}^{\infty}\primed
\frac{\exp\{-3it(k-j)(kj+m(k-j-m))\}}{kj(k-j)}
v_{k-j-m}(jv_{j})(mv_{m}),
\end{equation}
\begin{equation}\nonumber
R^{2}_{k}=\frac{i}{3}\sum_{j=-\infty}^{\infty}\primed
\sum_{\ell=-\infty}^{\infty}\primed
\frac{\exp\{-3itj(k(k-j)+\ell(j-\ell)\}}{kj(k-j)}v_{k-j}
\Bigg(j(v_{j-\ell})(\ell v_{\ell})\Bigg).
\end{equation}
We further rewrite $R^{2}_{k}$ by writing $j=j-\ell+\ell:$
\begin{multline}\nonumber
R^{k}_{2}=\frac{i}{3}\sum_{j=-\infty}^{\infty}\primed
\sum_{\ell=-\infty}^{\infty}\primed
\frac{\exp\{-3itj(k(k-j)+\ell(j-\ell))\}}{kj(k-j)}v_{k-j}
((j-\ell)v_{j-\ell})(\ell v_{\ell})
\\
+\frac{i}{3}\sum_{j=-\infty}^{\infty}\primed
\sum_{\ell=-\infty}^{\infty}\primed
\frac{\exp\{-3itj(k(k-j)+\ell(j-\ell))\}}{kj(k-j)}v_{k-j}v_{j-\ell}
(\ell^{2}v_{\ell}).
\end{multline}

Using the inequality $$\left|\frac{k^{2}}{kj(k-j)}\right|\leq 2,$$
and our previous arguments, we find the following bounds:
$$\|B\|_{H^{s+2}}\leq c\|v\|_{H^{s}}\|v\|_{H^{s+1}},$$
$$\|R^{1}\|_{H^{s+2}}\leq c\|v\|_{H^{s}}\|v\|_{H^{s+1}}^{2},$$
$$\|R^{2}\|_{H^{s+2}}\leq c\|v\|_{H^{s}}\left(\|v\|_{H^{s+1}}^{2}+\|v\|_{H^{s}}\|v\|_{H^{s+2}}\right).$$
If $\|v\|_{H^{s+2}}<1,$ then the right-hand sides can all be bounded by
$$c\|v\|_{H^{s}}\|v\|_{H^{s+2}}.$$
Following the arguments of the previous cases, then, we are able to conclude, for
sufficiently small $u_{0}\in H^{s+2}_{0},$
$$\|S_{D}(T)u_{0}\|_{H^{s+2}}\leq c\|u_{0}\|_{H^{s}}\|u_{0}\|_{H^{s+2}}.$$
This is estimate (\ref{smoothingGeneral}) with $p=\tilde{p}=2$ and $q=1,$ so we
see that the nonexistence result holds for the KdV equation:
\begin{cor}\label{mainSpecificResultKDV}
Let $0<T_{1}<T_{2}$ be given.  Let $p\in(1,2]$ be given.  Let $0<\delta<T_{2}-T_{1}$ be
given.  
Let $W_{p,\delta}\subseteq[T_{1},T_{2}]$ be as in Lemma \ref{linearLemma2},
with $A$ given by 
$$A=-\partial_{x}^{3}-\alpha\partial_{x},$$ 
with $|\alpha|\leq\frac{1}{2}.$
There exists $r_{1}>0$ such that 
for all $T\in W_{p,\delta},$ if
$u$ is a smooth,
non-constant time-periodic solution of (\ref{kdv1}) with temporal period $T,$ then 
$$\inf_{t\in[0,T]}\|u\|_{H^{6}}>r_{1}.$$

Let $W_{p}\subseteq[T_{1},T_{2}]$ be as in Corollary \ref{almostEveryCorollary}.  
Let $T\in W_{p}$ be given.  Then there exists
$r_{2}$ such that if $u$ is a smooth, non-constant time-periodic solution of (\ref{kdv1}) 
with temporal period $T,$ then
$$\inf_{t\in[0,T]}\|u\|_{H^{6}}>r_{2}.$$
\end{cor}

We note that the result of Corollary \ref{mainSpecificResultKDV} is given for solutions
in $H^{6}.$  This is because we take $s=4$ so that we work with classical solutions,
and we have $\tilde{p}=2.$  An application of Theorem \ref{generalTheorem} then gives
a restriction on the $H^{s+\tilde{p}}=H^{6}$ norm.

\section{Discussion}\label{section:discussion}

We have developed a theoretical framework for the demonstration of nonexistence of small
doubly periodic solutions for dispersive evolution equations.  The abstract theorem indicates that
nonexistence follows from the demonstration of dispersive smoothing estimates.  In particular cases,
we have demonstrated that the required dispersive smoothing estimates hold.  These results are
an analogue of scattering results for dispersive equations on the real line, since scattering 
implies the nonexistence of small-amplitude coherent structures.

Other work to be done includes treating additional specific examples, and possibly proving a
general theorem about when the weak smoothing property (\ref{smoothingGeneral}) holds.
For example, it should be investigated whether the necessary properties can be shown to
hold for other dispersive equations in one space dimension (like the Benjamin-Ono equation)
and in higher dimensions (such as Schr\"{o}dinger equations).  For equations with sufficiently
strong dispersion, the stronger smoothing estimate (corresponding to $\tilde{p}=0,$ in which the
Duhamel integral gains more than one derivative) will likely hold.
A class of equations with strong dispersion are fourth-order Schr\"{o}dinger equations
(see, e.g., \cite{fourthOrder2}, \cite{fourthOrder3}, \cite{fourthOrder}).   Such equations are of the form
$$i\psi_{t}+\Delta\psi+|\psi|^{2\sigma}\psi+\varepsilon\Delta^{2}\psi=0,$$
with $\sigma>0$ and $\varepsilon>0,$ and can arise by including higher-order corrections when
deriving a Schr\"{o}dinger equation from the Maxwell equations.
The above linear estimates, such as Lemma \ref{linearLemma2} and Corollary 
\ref{almostEveryCorollary}, are valid in one spatial dimension.  As pointed out to the authors by the 
referee of \cite{ambroseWright-CR}, in $n$ spatial dimensions, the result requires $p>n$ rather
than $p>1.$  Thus, in higher dimensions, one would expect to need to use the weaker smoothing 
property (corresponding to $\tilde{p}>0$) in order to follow the strategy of the
present paper.

\bibliographystyle{plain}
\bibliography{duhamel-smoothing}

\end{document}